\documentclass[10pt]{article}
\usepackage{amsmath}
\usepackage{amsfonts}
\usepackage{amssymb}
\usepackage{amsthm}
\usepackage{palatino}
\usepackage[margin=2.5cm, vmargin={1.5cm},includefoot]{geometry}

\title{Unipotent vector bundles and higher-order non-holomorphic
Eisenstein series}

\author{Jay Jorgenson\footnote{The first author was partially supported by NSF and PSC-CUNY grants.},
Cormac O'Sullivan\footnote{The second author was partially
supported by PSC-CUNY Research Award No. 67478-00 36 and PSC-CUNY Collaborative Award No. 80209-04 11.}}
\date{7 October 2006}

\begin{document}

\maketitle

%domains
\def\H{{\mathbf{H}}}
\def\F{{\mathbf F}}
\def\C{{\mathbf C}}
\def\R{{\mathbf R}}
\def\Z{{\mathbf Z}}
\def\Q{{\mathbf Q}}
\def\N{{\mathbf N}}
\def\B{{\mathbf B}}
%symbols
\def\G{{\Gamma}}
\def\GH{{\G \backslash \H}}
\def\g{{\gamma}}
\def\L{{\Lambda}}
\def\ee{{\varepsilon}}
\def\K{{\mathcal K}}
\def\Re{\mathrm{Re}}
\def\Im{\mathrm{Im}}
\def\PSL{\mathrm{PSL}}
\def\SL{\mathrm{SL}}
\def\Vol{\operatorname{Vol}}
\def\id{\operatorname{I}}

\def\sgn{\operatorname{sgn}}

%cusps
\def\ca{{\mathfrak a}}
\def\cb{{\mathfrak b}}
\def\cc{{\mathfrak c}}
\def\cd{{\mathfrak d}}
\def\ci{{\infty}}

%scaling matrices
\def\sa{{\sigma_\mathfrak a}}
\def\sb{{\sigma_\mathfrak b}}
\def\sc{{\sigma_\mathfrak c}}
\def\sd{{\sigma_\mathfrak d}}
\def\si{{\sigma_\infty}}

%modular symbol
\def\s#1#2{\langle \,#1 , #2 \,\rangle}

%theorems

\newtheorem{theorem}{Theorem}[section]
\newtheorem{lemma}[theorem]{Lemma}
\newtheorem{prop}[theorem]{Proposition}
\newtheorem{cor}[theorem]{Corollary}

\renewcommand{\labelenumi}{(\roman{enumi})}

\numberwithin{equation}{section}

\bibliographystyle{plain}

\begin{abstract}\noindent
Higher-order non-holomorphic Eisenstein series associated to a
Fuchsian group $\Gamma$ are defined by twisting the series
expansion for classical non-holomorphic Eisenstein series by
powers of modular symbols. Their functional identities  include
multiplicative and additive factors, making them distinct from
classical Eisenstein series.  In this article we prove the
meromorphic continuation of these series and establish their
functional equations which relate values at $s$ and $1-s$. In
addition, we construct high rank vector bundles $\cal V$ from
certain unipotent representations $\pi$ of $\Gamma$ and show that
higher-order non-holomorphic Eisenstein series can be viewed as
components of certain eigensections, $\mathbb E$, of $\cal V$.
With this viewpoint the functional identities of these
higher-order
 series are formally identical to the classical case.  Going
 further, we prove bounds for the Fourier coefficients of the
 higher-order non-holomorphic Eisenstein series.
\end{abstract}

\section{Introduction}

 Let $\Gamma$ be any Fuchsian group of the first kind which acts on
the hyperbolic upper half-space $\H$ such that the quotient
$\Gamma \backslash \H$ has finite volume yet is non-compact.
Following the notation from \cite{Iw}, let us fix representatives
for the finite number of $\G$-inequivalent cusps, label them
$\ca,\cb, \dots$, and use the scaling matrices $\sa,\sb, \dots$ to
give local coordinates near these cusps. Define the subgroup
$\G_\ca$ to be the  elements of $\G$ which fix the cusps
equivalent to $\ca$. Let $S_2(\G)$ be the space of holomorphic
cusp forms of weight 2 with respect to $\G$. For $f,g \in S_2(\G)$
and $m,n \geqslant 0$, we define, following \cite{Go2},
\cite{Go1}, and \cite{O'S1}, the higher-order non-holomorphic
Eisenstein series associated to the forms $f$ and $g$, character
$\chi$ and cusp $\ca$ by the series
\begin{equation}\label{eis}
    E_\ca^{m,n}(z,s;f,g,\chi)=\ \sum_{\g \in \G_\ca\backslash\G} \chi(\g) \s{\g}{f}^m
\overline{\s{\g}{g}}^n
\Im(\sa^{-1}\g z)^s
\end{equation}
where the modular symbols are defined by
\begin{equation}\label{mod}
    \s{\g}{f} = \int_{z_0}^{\g z_0} f(w) \, dw, \ \ \s{\g}{g} = \int_{z_0}^{\g z_0} g(w) \,
    dw\,;
\end{equation}
observe that the modular symbols are independent of $z_0 \in \H$
by Cauchy's theorem since $f$ and $g$ are holomorphic and weight
2. We will discuss the exact meaning of {\it order} in  section
\ref{hof}. Throughout this article, we will consider the forms
$f$, $g$ and character $\chi$ to be fixed, hence we will
abbreviate the notation and just write $E_\ca^{m,n}(z,s)$.

\vskip .10in As
will be shown, the series (\ref{eis}) converges
 absolutely for $\Re (s)
> 1$ and $z \in \H$.
One of the main results of this article is to prove the
meromorphic continuation of $E_\ca^{m,n}(z,s)$ to all $s$ in $\C$.
The continuation of (\ref{eis}) has already been addressed in
\cite{Pe} and \cite{PR} using perturbation theory.  Our new proof
continues the extension of Selberg's method \cite{Se} to
higher-order forms that appears in \cite{O'S1}, \cite{JO'S}. It
has the advantage of naturally giving strong bounds on both the
Fourier coefficients of $E_\ca^{m,n}(z,s)$ and its growth in $z$.
We require these bounds for our demonstration of the functional
equation, relating values of higher-order Eisenstein series in the
left and right $s$-planes; the determination of the functional
equation of (\ref{eis}) is one of the main new features of the
present article.  Finally, we formulate our results using the
language and notation of vector bundles, which shows that the
Eisenstein series (\ref{eis}) naturally can be viewed as
components of eigensections of certain vector bundles.

\vskip .10in The properties of $E_\ca^{m,n}$ are developed here
with two new applications in mind.  First, it is well known that
the Kronecker limit formula for the classical Eisenstein series
$E_\ca^{0,0}$ (usually written as just $E_\ca$) produces Dedekind
sums. In our forthcoming work \cite{JO'S2}, we find that
$E_\ca^{1,1}$ yields new types of Dedekind sums with interesting
arithmetic. No such sums were found for $E_\ca^{1,0}$ though the
explicit Kronecker limit formula for that case was revealed in
\cite{JO'S} which included new families of $L$-functions
associated to the Fourier coefficients of modular forms and
classical Kronecker limit functions. Second, in current work of
the second author and Imamo\={g}lu, a natural inner product for
the space of higher-order weight $k$ forms is found. The usual
Petersson inner product may be unfolded to an integral over a
vertical strip. For the higher-order inner product the strip is
refolded and detailed knowledge of $E_\ca^{m,n}(z,s)$ is required.

\vskip .10in We view the ideas developed in this article as
establishing fundamental results in the theory of higher-order
non-holomorphic Eisenstein series.  Numerous other authors have
encountered higher-order modular forms in their research, and we
mention briefly three applications that have appeared. Firstly, in
\cite{Go2}, \cite{Go1} Goldfeld defined the series $E_\ca^{m,n}$
for the first time and shows their utility in studying the
distribution of modular symbols. Detailed results of this type are
required to provide an approach to the $ABC$-conjecture via
Goldfeld's conjectures on periods and modular symbols, as
described in \cite{Go4}, \cite{Go3}.  Secondly, in \cite{PR},
\cite{Ri}. Petridis and Risager prove that, under a certain
ordering and normalization, modular symbols follow a Gaussian
distribution. Finally, in \cite{DO'S} it is shown that applying
the Hecke operators to the residues of poles of $E_\ca^{1,0}$ on
the critical line $\Re(s)=1/2$ yields an identity that can be used
to construct an orthonormal basis for a certain space of Maass
forms.

\section{Statement of results}
Let $\mathcal E^{m,n} = (E^{m,n}_{\ca})$ be the $r \times 1$
column vector of Eisenstein series, where we index over all $r$
inequivalent cusps $\ca$ of $\G$. Let $\mathcal I_{m,n}$ denote
all pairs $(i_1,i_2)$ for $0 \leqslant i_1 \leqslant m$ and  $0
\leqslant i_2 \leqslant n$, and with the ordering which has
$(i_1,i_2) < (j_1,j_2)$ when $i_1 < j_1$ or when $i_1 = j_1$ and
$i_2 < j_2$.
 Let $N=r(m+1)(n+1)$ and define the $N \times 1$
column vector
\begin{equation}\label{eee}
\mathbb E^{m,n}= {}^t\left( {}^t \mathcal E^{m-i_1,n-i_2}
\right)_{(i_1,i_2) \in \mathcal I_{m,n}}.
\end{equation}
For example, for a group $\G$ with $r=3$ inequivalent cusps, we
have when $(m,n) = (2,1)$,
\begin{eqnarray*}
\mathbb E^{2,1} & = & {}^t\left( {}^t \mathcal E^{2,1}, {}^t
\mathcal E^{2,0}, {}^t
\mathcal E^{1,1}, {}^t \mathcal E^{1,0}, {}^t \mathcal E^{0,1}, {}^t \mathcal E^{0,0}  \right),\\
& = & {}^t\left(  E^{2,1}_\ca, E^{2,1}_\cb, E^{2,1}_\cc,
E^{2,0}_\ca, E^{2,0}_\cb, E^{2,0}_\cc, \cdots , E^{0,0}_\ca,
E^{0,0}_\cb, E^{0,0}_\cc \right).
\end{eqnarray*}
In this notation, we now state our first main result.

\begin{theorem} \label{thm1} There exists an explicit representation
$\pi^{m,n}$ of $\Gamma$ into $\text{\rm Uni}(N, \mathbf C)$, the
space of $N \times N$ unipotent matrices, such that for all
$\gamma \in \Gamma$ and $\Re (s) > 1$, we have
\begin{equation}\label{rep}
\mathbb E^{m,n}(\gamma z,s) = \pi^{m,n}(\gamma)\mathbb
E^{m,n}(z,s).
\end{equation}
Furthermore
\begin{equation}\label{del}
\Delta \mathbb E^{m,n}(z,s) = s(1-s) \mathbb E^{m,n}(z,s),
\end{equation}
 where $\Delta$ denotes the usual hyperbolic
Laplacian $-4y^2(d/dz)(d/d\overline{z})$.
\end{theorem}

\vskip .10in One can restate Theorem \ref{thm1} as follows.  There
exists a unipotent representation $\pi^{m,n}$ of $\Gamma$ and
associated vector bundle ${\cal V}_{m,n}$ over $\Gamma \backslash
\H$ such that for $\Re (s) > 1$, the vector $\mathbb E^{m,n}$ is a
$C^{\infty}$ section of ${\cal V}_{m,n}$; in addition, $\mathbb
E^{m,n}$ is an eigenvector under the action by the Laplacian
$\Delta$ with eigenvalue $s(1-s)$.  If $(m,n) = (0,0)$, then
Theorem \ref{thm1} reduces to  well-known properties of the
classical non-holomorphic Eisenstein series.  For all other pairs
$(m,n)$, Theorem \ref{thm1} encodes the multiplicative and
additive behavior of the higher-order non-holomorphic Eisenstein
series when acted upon by the group $\Gamma$ and the Laplacian
$\Delta$.

\vskip .10in Theorem \ref{thm1} shows that $\mathbb E$ is a
vector-valued automorphic form in the sense of Knopp and Mason
\cite{KM1}, \cite{KM2}. Non-holomorphic Eisenstein series
associated to monomial representations (matrices with one non-zero
entry in each row and column) are studied in \cite{Ta}. These
series arise naturally when attempting to extend the
Rankin-Selberg method, and Taylor in \cite{Ta} has shown the
meromorphic continuation of such series.

\vskip .10in Our next result determines the meromorphic
continuation and functional equation for the vector $\mathbb
E^{m,n}$. Recall that $\s{\g}{f}=\s{\g}{g}=0$ for $\g \in \G_\cb$,
since $f$ and $g$ are cusp forms, so then, for any cusps $\ca$ and
$\cb$, we have that $E_{\ca}^{m,n}(\sb(z+1) ,s)=
E_{\ca}^{m,n}(\sb(z) ,s)$. We conclude that the Eisenstein series
$E_{\ca}^{m,n}$ admits a Fourier expansion at the cusp $\cb$ which
can be shown to be
\begin{equation}\label{e2}
E_{\ca}^{m,n}(\sb z,s)=\delta_{0,0}^{m,n}\cdot\delta_{\ca \cb}
y^s+ \phi_{\ca \cb}^{m,n}(s)y^{1-s}+ \sum_{k\not=0}\phi_{\ca
\cb}^{m,n}(k,s)W_s(kz),
\end{equation}
where $z=x+iy$,
$$
W_s(kz)=2\sqrt{|k|y}K_{s-1/2}(2\pi |k| y)e^{2\pi
i k x}
$$ and $K_{s-1/2}$ is the $K$-Bessel function.  The factor
$\delta_{0,0}^{m,n}$ is 1 if $(m,n)=(0,0)$ and 0 otherwise. Also
$\delta_{\ca \cb}$ is 1 if $\ca = \cb$ and 0 otherwise.   Equation
(\ref{e2}) reveals an important property of $E^{m,n}$: For $(m,n)
\neq (0,0)$ and $\Re(s)>1$, $E^{m,n}$ has polynomial decay at
every cusp. In essence, this holds because of the definition
(\ref{eis}) and since $\s{\id_2}{f}=\s{\id_2}{g}=0$ where we use
$\id_k$ to denote the $k \times k$ identity matrix.

\vskip .10in For $(m,n)=(0,0)$, the expansion
(\ref{e2}) is proved in  \cite[(3.20)]{Iw} and that proof generalizes to
 the case $(m,n)=(1,0)$ in \cite{O'S1} and easily
 to all other $(m,n)$.  Furthermore, the  coefficients $\phi_{\ca
\cb}^{m,n}(k,s)$ can be identified in terms of Kloosterman sums
twisted by modular symbols.

\vskip .10in For fixed $i$ and $j$, let us assemble the functions
$\phi_{\ca \cb}^{i,j}(s)$ in the constant term (in $x$) of
(\ref{e2}) into an $r \times r$ matrix
$$
\Phi^{i,j}(s) =\left( \phi_{\ca \cb}^{i,j}(s)\right).
$$
We now construct $\Phi_\mathbb E^{m,n}$, an $N \times N$ matrix
from the $r \times r$ matrices $\Phi^{i,j}$ by setting
$$
\Phi_\mathbb E^{m,n}=\left( \binom{m-i_1}{j_1-i_1}
\binom{n-i_2}{j_2-i_2} \Phi^{j_1-i_1,
j_2-i_2}\right)_{\bigl((i_1,i_2),(j_1,j_2) \bigr) \in \mathcal
I_{m,n} \times \mathcal I_{m,n}}
$$
where we understand that $\Phi^{u,v}=0$ if $u$ or $v$ is negative.
Note that we are following standard conventions for binomial
coefficients which set $\binom{0}{0}=1$ and $\binom{n}{k} =0$ if
$k<0$ or $k>n$.

\vskip .10in
\begin{theorem} \label{thm2} The vector of Eisenstein series
$\mathbb E^{m,n}(z,s)$ admits a meromorphic continuation to all $s
\in \mathbf C$ and satisfies the functional equation
\begin{equation}\label{funeq}
\mathbb E^{m,n}(z,1-s)=\Phi_{\mathbb E}^{m,n}(1-s)\mathbb
E^{m,n}(z,s).
\end{equation}
In addition, the matrix $\Phi_\mathbb E^{m,n}$ satisfies the
identity
\begin{equation}\label{funeq2}
\Phi_{\mathbb E}^{m,n}(1-s)\Phi_{\mathbb E}^{m,n}(s)=\id_N\,.
\end{equation}
\end{theorem}

\vskip .10in In the case when $(m,n)=(0,0)$, Theorem \ref{thm2} is
a classical result of Selberg \cite{Se}. When $(m,n)=(1,0)$, the
proof of Theorem \ref{thm2} appears in \cite{O'S1} where Selberg's
original methods in the case $(m,n)=(0,0)$ were modified, taking
into account that $\mathcal E^{1,0}$ is no longer
$\Gamma$-invariant. However, the notation in \cite{O'S1} does not
reveal the structure implied by the notation in Theorem \ref{thm1}
or Theorem \ref{thm2}, which is one important part of the point of
view taken here.  Without the notation of Theorem \ref{thm2}, the
functional equations (\ref{funeq}) and (\ref{funeq2}) can be shown
to be equivalent to the formulas
\begin{eqnarray}\label{funeq3}
 \sum_{0 \leqslant i \leqslant m \atop
0 \leqslant j \leqslant n} \binom{m}{i} \binom{n}{j}
\Phi^{i,j}(1-s) \mathcal E^{m-i,n-j}(z,s) & = & \mathcal E^{m,n}(z,1-s), \\
\sum_{0 \leqslant i \leqslant m \atop
0 \leqslant j \leqslant n} \binom{m}{i} \binom{n}{j}
\Phi^{i,j}(1-s) \Phi^{m-i,n-j}(s) & = & \begin{cases}
\id_r & \text{if \ \ $m=n=0$,} \\
0 & \text{otherwise.}
\end{cases} \label{phieq}
\end{eqnarray}
These results should be compared with those of Risager in
\cite[Chapter 4]{Ri} where a functional equation is found with
matrix entries involving perturbation coefficients. These
coefficients are also expressed there in terms of certain
Dirichlet series.

\vskip .10in The next result is an example of the growth
properties of $E_\ca^{m,n}$ which we prove in section
\ref{continuation}. We shall use the notation $s=\sigma + it$ throughout the paper.

\vskip .10in
\begin{theorem} \label{thm3} For each compact set ${\mathbf S} \subset
\C$ and all integers $m,n \geqslant 0$, there exists a
holomorphic function $\xi^{m,n}_{\mathbf S} (s)$ such that for  all $s$ in ${\mathbf
S}$ and all $k \neq 0$, we have
\begin{equation}
\xi^{m,n}_{\mathbf S}(s)\phi^{m,n}_{\ca \cb}(k,s)  \ll  (\log^{m+n} |k| +1)(|k|^\sigma + |k|^{1-\sigma}). \label{phi71}
\end{equation}
The implied constant depends solely on $f, g$ and $\G$.
\end{theorem}

\vskip .10in Finally, we note that the most general form of our
theory of higher-order non-holomorphic Eisenstein series includes
more that two modular symbols. As shown in \cite{GO'S}, for
example, the space of all homomorphisms from $\G$ to $\C$ that are
zero on parabolic elements is $2g$ dimensional, where $g$ is the
genus of $\GH$, and is generated by modular symbols and their
conjugates. If $h_i$ is a basis of such homomorphisms, then we may
consider the general series
\begin{equation}\label{geneeis}
E_\ca^{m_1,m_2, \cdots, m_{2g}}(z,s;h_1,h_2, \cdots,
h_{2g},\chi)=\ \sum_{\g \in \G_\ca\backslash\G} \chi(\g)
h_1(\g)^{m_1} h_2(\g)^{m_2} \cdots h_{2g}(\g)^{m_{2g}}
\Im(\sa^{-1}\g z)^s.
\end{equation}
As the reader may verify, the methods in this article apply with
very few changes to the series (\ref{geneeis}). Homomorphisms that
are not zero on all parabolic elements can also be used in
(\ref{geneeis}); see \cite{GO'S} for some initial results on these
series.

\section{Some preliminary results}

\subsection{Higher-order forms} \label{hof}

Let us set $\mathcal A^0(\G)=\{0\}$ and, for $n \geqslant 1$,
define $\mathcal A^n(\G)$ to be the space of smooth functions
$\psi: \H \rightarrow \C$ such that
$$
\psi(\gamma z) -\psi(z) \in \mathcal A^{n-1}(\G) \text{ \ \ \ for
all \ }\gamma \in \G.
$$
If $n=1$, then $\mathcal A^1(\G)$ is the well-known space of
automorphic functions.    We call the elements of $\mathcal
A^n(\G)$ $n$th-order automorphic forms.  In this general setting,
$n$th-order automorphic forms were first described in \cite{KZ}
and \cite{CDO}.

\vskip .10in More completely, let $j(\g,z)$ be the usual
automorphy factor associated to automorphic forms.  For $\g \in
\G$, character $\chi$ of $\G$, and integer $k \in \Z$,  the slash
operator is defined on a function $f$ by
$$
(f|_{k,\chi} \g)(z)=\frac{f(\g z)}{j(\g,z)^k \chi(\g)}.
$$
Recursively we define $\mathcal A^n_k(\G, \chi)$ to be the space of
all smooth functions such that
$$
\psi |_{k,\chi} \g -\psi \in \mathcal A^{n-1}_k(\G,\chi) \text{ \ \ \
for all \ }\g \in \G,
$$
where, as before, we set $\mathcal A^0_k(\G, \chi)=\{0\}$. In
words, $\mathcal A^n_k(\G, \chi)$ is the space of weight $k$,
$n$th-order automorphic forms which are twisted by the character
$\chi$.  Observe that if we set $\xi_\g = \psi |_{k,\chi} \g
-\psi$ then $\xi_\g$, as the notation indicates, may depend on
$\g$, though not arbitrarily:  we necessarily have that $\xi_{\g
\tau} = \xi_\g |_{k,\chi} \tau + \xi_\tau$ for all $\g, \tau \in
\G$.

\begin{lemma} \label{subsp} For any $0 \leqslant n \leqslant m$ we
have $\mathcal A^n_k(\G,\chi) \subseteq \mathcal A^m_k(\G,\chi)$.
\end{lemma}
\begin{proof}
We use induction on $n$.
Trivially, we have that $\mathcal A^0_k(\G,\chi) \subseteq \mathcal
A^m_k(\G,\chi)$ for any $m \geqslant 0$. This is the base case.
For the inductive step:
\begin{eqnarray*}
\psi \in \mathcal A^n_k(\G,\chi) & \implies & \psi |_{k,\chi} \g -\psi \in \mathcal A^{n-1}_k(\G,\chi)\\
& \implies & \psi |_{k,\chi} \g -\psi \in \mathcal A^{m-1}_k(\G,\chi), \ \ m\geqslant n \,\,\,
{{\textrm{(by induction)}}}\\
& \implies & \psi  \in \mathcal A^{m}_k(\G,\chi).
\end{eqnarray*}
This proves the claim.
\end{proof}

\begin{prop} \label{prod} For any $n,m \geqslant 1$,
if $f \in \mathcal A^n_k(\G,\chi_1)$ and  $g \in \mathcal A^m_l(\G,\chi_2)$
then
$$
f\cdot g \in \mathcal A^{n+m-1}_{k+l}(\G,\chi_1 \chi_2).
$$
\end{prop}
\begin{proof}
 We use
induction on $n+m$. Observe that
$$
\left(f\cdot g \right) |_{k+l,\chi_1\chi_2 } \g = (f |_{k,\chi_1} \g)( g |_{l, \chi_2}
\g) =(f+f_1)(g+g_1),
$$
for some $f_1 \in \mathcal A^{n-1}_k(\G,\chi_1)$ and $g_1 \in \mathcal
A^{m-1}_l(\G,\chi_2)$. Hence
$$
\left(f\cdot g \right) |_{k+l,\chi_1\chi_2 } \g - f\cdot g=f_1 \cdot g+ f\cdot g_1 + f_1 \cdot g_1.
$$
By induction and Lemma \ref{subsp}, we have that both $f_1 \cdot g$ and
$f\cdot g_1$ are elements of $\mathcal A^{n+m-2}_{k+l}(\G,\chi_1\chi_2)$,
and $ f_1\cdot g_1 \in \mathcal A^{n+m-3}_{k+l}(\G,\chi_1\chi_2) \subseteq
\mathcal A^{n+m-2}_{k+l}(\G,\chi_1\chi_2)$. Therefore, $f\cdot g
|_{k+l,\chi_1\chi_2 } \g - f\cdot g \in \mathcal A^{n+m-2}_{k+l}(\G,\chi_1\chi_2)$ which
implies that $f\cdot g \in \mathcal
A^{n+m-1}_{k+l}(\G,\chi_1\chi_2)$ as required. This proves the
induction step.

\vskip .10in It remains to prove the first step in the induction,
namely when $(n,m) = (1,1)$.  This follows immediately from
 the observation that the product of two classical
automorphic forms of weight $k$ and $l$ is an automorphic form of
weight $k+l$.   With all this, the proof of the proposition is
complete.
\end{proof}

\vskip .10in {\bf Remark:} We can interpret Lemma \ref{subsp} and
Proposition \ref{prod} as proving that the space of all
higher-order automorphic forms has a graded ring structure. In
\cite{CD} it is shown that the theory of  Rankin-Cohen brackets
applies and gives the set of higher-order forms a {\it canonical
RC structure}.  It is natural to ask if, under either structure,
the graded ring of higher-order holomorphic automorphic forms is
finitely generated.

\subsection{Higher-order Eisenstein series}

\vskip .10in We will show that $E_\ca^{m,n} \in \mathcal
A^{m+n+1}_0(\G,\chi)$ in Lemma \ref{eina}. In this section, we
establish some more basic properties.

\vskip .10in To see that the terms in (\ref{eis}) are well defined
we note that $\s{\g}{f} = 0$ whenever $\g \in \G_\ca$.  Therefore,
the modular symbols are well-defined functions on the cosets
$\G_{\ca} \backslash \G$.  Similarly, we also require that
$\chi(\g)=1$ for all $\g \in \G_\ca$ and for convergence the
character must satisfy $|\chi|=1$.

\vskip .10in The size of an automorphic form in cuspidal zones is
of great importance.  Following \cite{Iw}, we will measure the
growth of automorphic forms using the  {\it invariant height} function
defined by
$$
y_\G(z)=\max_\ca \max_{\g \in \G}(\Im( \sa^{-1}\g z)).
$$
Suppose that $\psi$ is a continuous function on $\H$ such that
$|\psi|$ is
 $\G$-invariant, and let $B$ denote a continuous positive function on
$[0,\infty)$. Then it is easy to verify that the following
conditions are equivalent:
\begin{enumerate}
\item $\psi(\sa z) \ll B(y)$ as $y \to \infty$ for all cusps $\ca$,
\item $\psi(z) \ll B(y_\G(z)).$
\end{enumerate}

 More generally, for continuous functions $\xi$ on $\H$ which are not
necessarily automorphic, it is more convenient to fix a
fundamental domain, $\F$, and examine their growth on $\F$. Let
$\mathcal P_Y \subset \H$ denote the strip with $|x| \leqslant
1/2$ and $y \geqslant Y$. We choose $\F$ so that it contains the
cuspidal zones $\sa \mathcal P_Y$ for all $\ca$ and $Y$ large
enough (see \cite[Section 2.2]{Iw}). For $z \in \F$ we define
$$
y_\F(z)=\max_\ca (\Im( \sa^{-1} z))
$$
to be the {\it domain height} function; with this definition, we
assert that the following statements are equivalent
\begin{enumerate}
\item $\xi(\sa z) \ll B(y)$ as $y \to \infty$ for all cusps $\ca$ and all $|x| \leqslant 1/2$,
\item $\xi(z) \ll B(y_\F(z))$ for all $z \in \F$.
\end{enumerate}
Also we observe the relation
$$
y_\F(z) \leqslant y_\G(z).
$$
Further required properties of the invariant height   are proved in Appendix \ref{app}.

\vskip .10in
\begin{prop} \label{prop23}
For $m,n \geqslant 0$ and $(m,n) \neq (0,0)$ the series (\ref{eis}
defining $E_\ca^{m,n}( z,s)$ converges to a smooth function of $z
\in \H$ and holomorphic function of $s \in \C$ provided $\Re (s) >
1$.  The convergence is absolute for any $s \in \C$ with $\Re(s) >
1$ and uniform in the region $\Re(s) \geqslant 1+\delta$ for any
$\delta>0$. Furthermore, for $\Re (s) = \sigma
> 1$, the function $E_\ca^{m,n}( z,s)$ satisfies the growth
condition
\begin{equation}\label{four}
E_\ca^{m,n}( z,s) \ll y_\F(z)^{1-\sigma+\varepsilon}
\end{equation}
 for any $\varepsilon>0$ and all $z \in \F$. The implied constant depends only on
 $\varepsilon, \sigma, m,n,f,g$ and $\G$.
\end{prop}

\vskip .10in
\begin{proof}
Set
\begin{equation}\label{fg}
F_\ca(z)=\int_\ca^z f(w)\,dw, \ \ G_\ca(z)=\int_\ca^z g(w)\,dw.
\end{equation}
First we note that, for any two (possibly equal) cusps $\ca$ and $\cb$ and for all $z \in \H$,
\begin{equation}\label{fab}
F_\ca(\sb z) \ll |\log y| +1
\end{equation}
 with the same bound for $G_\ca$, where the
implicit constant in (\ref{fab}) depends solely on $f$ and $\G$.
To prove (\ref{fab}), observe that the integral of $f(z)$ from $z$
to $z+1$ is zero, and that $f(z) \ll 1/y$
(see \cite[Lemma 3]{DKMO}). It also follows from (\ref{fab}) that
\begin{equation}\label{fab2}
F_\ca(z) \ll \log(y_\F(z) +e).
\end{equation}

\vskip 0.10in
We shall need some elementary inequalities in the sequel:
\begin{eqnarray}
|\log y| & < & (y^\epsilon+y^{-\epsilon})/\epsilon, \label{log}\\
|\log y|+1 & < & 3 \log (y+1/y), \label{log2}\\
(y+1)^r & < & 2^r (y^r +1), \label{yr}\\
(y+1/y)^r & < & 2^r (y^r +y^{-r})\label{yr2}
\end{eqnarray}
for all $y$, $\epsilon$, $r >0$.

\vskip 0.10in
From (\ref{fab}) and (\ref{log}), we deduce that
\begin{eqnarray}
\s{\g}{f} & = & F_\ca(\g z) - F_\ca( z) \nonumber\\
& \ll & |\log \Im(\sb^{-1} \g z)| + |\log \Im(\sb^{-1} z)|+1 \nonumber\\
&  \ll &  \Im(\sb^{-1} \g z)^\varepsilon + \Im(\sb^{-1} \g z)^{-\varepsilon}
 +  \Im(\sb^{-1} z)^\varepsilon +\Im(\sb^{-1} z)^{-\varepsilon}. \nonumber
\end{eqnarray}
Hence, for any $\varepsilon >0$, any $z \in \H$ and any cusp $\ca$,
\begin{equation} \label{modsym}
    \s{\g}{f}^m\s{\g}{g}^n \ll   \Im(\sa^{-1} \g z)^\varepsilon + \Im(\sa^{-1} \g z)^{-\varepsilon}
 +  \Im(\sa^{-1} z)^\varepsilon +\Im(\sa^{-1} z)^{-\varepsilon}
\end{equation}
with an implied constant depending on $f,g,\G,m,n$ and $\varepsilon$. Also for any cusp $\ca$ and any $z \in \F$
\begin{equation}\label{bd1}
\Im( \sa^{-1} z) \leqslant y_\F(z)
\end{equation}
by definition. It is well-known that the classical non-holomorphic
Eisenstein series, meaning (\ref{eis}) for $(m,n) = (0,0)$,
converges for $\Re(s) >1$ (uniformly for $\Re(s) \geqslant
1+\delta$) and satisfies the bound
\begin{equation}\label{bd2}
E_\ca(z,s) \ll y_\G(z)^\sigma;
\end{equation}
see, for example, \cite[Corollary 3.5]{Iw}. Going further, we have
that
\begin{equation}\label{eisbd}
E_\ca(z,s) - \Im (\sa^{-1} z)^s  \ll y_\F(z)^{1-\sigma},
\end{equation}
for $z \in \F$; for the sake of completeness, we will prove
(\ref{eisbd}) in detail in Lemma \ref{eminus}. Since $\s{\id_2}{f}
= \s{\id_2}{g} =0$, the term in (\ref{eis}) corresponding to the
identity coset
 vanishes for $(m,n) \neq (0,0)$.  With all this, we
can substitute the bound (\ref{modsym}) into the definition
(\ref{eis}) and use (\ref{bd1}) and (\ref{eisbd}) to complete the
proof of the proposition.
\end{proof}

\section{Unipotent higher-order vector bundles}

\subsection{Automorphy factors}

 For any $f \in S_{2}(\Gamma)$ it is an easy consequence of
elementary complex analysis that, for all
$\gamma_{1}, \gamma_{2} \in \Gamma$,
\begin{equation}\label{symb}
\s{\g_1 \g_2}{f}= \s{\g_1}{f}+\s{\g_2}{f}  .
\end{equation}
This is the basis for the following
computation.

\vskip .10in
\begin{lemma}\label{transform} For any $\tau \in \Gamma$, let
$$
S_{i,j}(\tau)=\overline{ \chi}(\tau) \bigl(-\s{\tau}{f}\bigr)^{i}
\bigl(-\overline{\s{\tau}{g}}\bigr)^{j}.
$$
For $\Re (s) > 1$, the Eisenstein
series $E^{m,n}_{\ca}$ obeys the transformation rule
\begin{equation}\label{emn}
E^{m,n}_{\ca}(\tau z,s)=\sum_{i=0}^m\sum_{j=0}^n
\binom{m}{i}\binom{n}{j}S_{m-i,n-j}(\tau)E^{i,j}_{\ca}(z,s).
\end{equation}
\end{lemma}
\begin{proof}
Combining (\ref{eis}) and (\ref{symb}), we get
\begin{eqnarray*}
 E_\ca^{m,n}(\tau z,s) &=& \sum_{\g \in
\G_\ca\backslash\G}
\chi(\g) \s{\g}{f}^m \overline{\s{\g}{g}}^n
\Im(\sa^{-1}\g \tau z)^s
\\
& = &\sum_{\g \in
\G_\ca\backslash\G}
\chi(\g\tau^{-1})\s{\g\tau^{-1}}{f}^m
\overline{\s{\g\tau^{-1}}{g}}^n \Im(\sa^{-1}\g  z)^s
\\
& = &\sum_{\g \in
\G_\ca\backslash\G}
\chi(\g) \overline{ \chi}(\tau)\bigl(\s{\g}{f} - \s{\tau}{f} \bigr)^m \cdot
\bigl(\overline{\s{\g}{g}}- \overline{\s{\tau}{g}}\bigr)^n \Im(\sa^{-1}\g
z)^s.
\end{eqnarray*}
The result now follows by expanding the terms
$
(\s{\g}{f} - \s{\tau}{f} )^m $ and $
(\overline{\s{\g}{g}}- \overline{\s{\tau}{g}})^n$.
\end{proof}

\begin{lemma} \label{eina} For $m,n \geqslant 0$ we have $E_\ca^{m,n} \in \mathcal
A^{m+n+1}_0(\G,\chi)$.
\end{lemma}
\begin{proof}
Using (\ref{emn}) we observe that $E_\ca^{0,0}|_{0, \chi}\g -
E_\ca^{0,0}=0$ and $E_\ca^{0,0} \in \mathcal A^{1}_0(\G,\chi)$.
Now, using induction on $m+n$, we have by (\ref{emn}) that
$E_\ca^{m,n}|_{0, \chi}\g - E_\ca^{m,n}$ is a linear combination
of lower order Eisenstein series that are in $\mathcal
A^{m+n}_0(\G,\chi)$. The lemma follows.
\end{proof}

\vskip .10in We now complete the proof of Theorem \ref{thm1}.

\vskip .10in
\begin{proof}
Directly from (\ref{emn}) we obtain the relation
\begin{equation}\label{vectorident}
{\mathcal E}^{m,n} (\g z,s)=\sum_{i=0}^m\sum_{j=0}^n
\binom{m}{i}\binom{n}{j}S_{m-i,n-j}(\g){\mathcal E}^{i,j}(z,s)
\end{equation}
 for the $r \times 1$ vector of Eisenstein series
${\mathcal E}^{m,n}$. Construct the $N \times N$ matrix
\begin{equation}\label{constr}
\pi^{m,n}(\g) =\left( \binom{m-i_1}{m-j_1}\binom{n-i_2}{n-j_2} S_{j_1-i_1, j_2-i_2}(\g)
\id_r \right)_{\bigl((i_1,i_2),(j_1,j_2) \bigr) \in \mathcal I_{m,n} \times \mathcal I_{m,n}}.
\end{equation}
If we replace $(m,n)$ in (\ref{vectorident}) by $(u,v)$, say, then
it is elementary to verify that (\ref{vectorident}) can be
combined for all $(u,v) \leqslant (m,n)$ in the single equation
\begin{equation}\label{e7}
\mathbb E^{m,n}(\gamma z,s) = \pi^{m,n}(\gamma)\mathbb
E^{m,n}(z,s)\,,
\end{equation}
which holds for all $\g \in \Gamma$, proving
(\ref{rep}).  As for (\ref{del}), this relation following directly
from (\ref{eis}) when differentiating term by term, using that
$\Delta y^s = s(1-s)y^s$ and the fact that $\Delta$ is $\SL_2(\R)$
invariant.

\vskip .10in The matrix $\pi^{m,n}(\tau)$ is easily seen to be
upper triangular with ones on the diagonal, hence it is unipotent.
Upon replacing $z$ by $\g z$ and repeating (\ref{e7}), we obtain
the relation $\pi^{m,n}(\tau \g) = \pi^{m,n}(\tau)\pi(\g)$.
Therefore $\pi^{m,n}$ is a representation of $\Gamma$ into the set
of $N \times N$ unipotent matrices with $N=r(m+1)(n+1)$.
Consequently, the representation $\pi^{m,n}$ can be used to define
a rank $N$ unipotent vector bundle ${\cal V}^{m,n}$ over $\GH$,
and (\ref{rep}) is simply the statement that $\mathbb
E^{m,n}(z,s)$ is a section of ${\cal V}^{m,n}$.

\vskip .10in With all this, the proof of Theorem \ref{thm1} is
complete.
\end{proof}

\subsection{An example}
 A simple case
of Theorem \ref{thm1} occurs when $(m,n)=(2,0)$, for an arbitrary
number of cusps $r$.  In this instance, we see from (\ref{constr}) that
\begin{equation} \label{eg1}
\pi^{2,0}(\g) = \overline{\chi}(\g) \left( \begin{matrix}
\id_r & -2\s{\g}{f} \cdot \id_r & \s{\g}{f}^2 \cdot \id_r \\  &
\id_r & -\s{\g}{f} \cdot \id_r \\
 & &  \id_r \end{matrix}\right)
\end{equation}
and (\ref{e7}) says $\mathbb E^{2,0}(\gamma z,s) =
\pi^{2,0}(\gamma)\mathbb E^{2,0}(z,s)$ which encodes the three
identities

\begin{eqnarray*}
\mathcal E^{2,0}(\g z,s) & = & \overline{\chi}(\g) \left(\mathcal
E^{2,0}(z,s) -2\s{\g}{f} \mathcal
E^{1,0}(z,s) + \s{\g}{f}^2 \mathcal
E^{0,0}(z,s)\right),\\
\mathcal E^{1,0}(\g z,s) & = & \overline{\chi}(\g) \left(\mathcal
E^{1,0}(z,s) - \s{\g}{f} \mathcal
E^{0,0}(z,s)\right),\\
\mathcal E^{0,0}(\g z,s) & = & \overline{\chi}(\g) \mathcal
E^{0,0}(z,s)\, .
\end{eqnarray*}

\section{Constructing an integral equation for $E_\ca^{m,n}$}

\subsection{Bounds for $E_\ca^{m,n}$ on $\H$}

 We need bounds for $E_\ca^{m,n}$ on the upper half plane $\H$.
 The fact that it is not $\G$-invariant makes this hard to do directly,
 so we introduce a related automorphic series.

\vskip .10in For non-negative integers $m$ and $n$, we define,
using the notation in (\ref{fg}), the function
\begin{equation}\label{qmn}
Q_\ca^{m,n}(z,s) = \sum_{\g \in \G_\ca \backslash \G} \chi(\g)
F_\ca(\g z)^m \overline{G_\ca(\g z)}^n \Im(\sa^{-1} \g z)^s
\end{equation}
which will naturally be automorphic wherever it converges. Observe
that when $m=n=0$ the series (\ref{qmn}) is just the usual
Eisenstein series $E_\ca(z,s)$. To see the connection between
$E_\ca^{m,n}$ and $Q_\ca^{m,n}$ use the identity $F_\ca(\g z)
-F_\ca(z) = \langle \g, f \rangle$ and the binomial theorem to
obtain
\begin{eqnarray}\label{qe}
Q_\ca^{m,n}(z,s) & = & \sum_{i,j} \binom{m}{i} \binom{n}{j}
F_\ca(z)^{m-i} \overline{G_\ca(z)}^{n-j} E_\ca^{i,j}(z,s), \\
E_\ca^{m,n}(z,s) & = & \sum_{i,j} \binom{m}{i} \binom{n}{j}
\left(-F_\ca(z)\right)^{m-i} \left(-\overline{G_\ca(z)}\right)^{n-j}Q_\ca^{i,j}(z,s). \label{qe2}
\end{eqnarray}

Since $E^{i,j}_\ca(z,s)$ is convergent for $\Re(s)>1$, so is
$Q^{m,n}_\ca(z,s)$ by (\ref{qe}). Using (\ref{four}), (\ref{fab2})
and (\ref{bd2}) to bound the right side of (\ref{qe}), we get
\begin{equation}\label{qbound}
Q_\ca^{m,n}(z,s)  \ll  \log^{m+n}(y_\F(z)+e) \cdot y_\F(z)^\sigma
\end{equation}
where the  implied constant depends on $\epsilon$, $\sigma$, $m$,
$n$, $f$, $g$ and $\G$ alone and $z$ is restricted to $\F$. Since
$Q_\ca^{m,n}$ is automorphic, and the right side of (\ref{qbound})
is an increasing function of $y_\F$, it follows that, for all $m,n
\geqslant 0$ and $z \in \H$,
\begin{equation}\label{qmn2}
Q_\ca^{m,n}(z,s) \ll \log^{m+n}(y_\G(z)+e) \cdot  y_\G(z)^\sigma.
\end{equation}
We need the next result, which is a slightly stronger and more explicit version of \cite[Lemma 7.1]{JO'S}.

\begin{lemma} \label{bnd71} Suppose $D$ is a
function on $\H$ and $B$ is an increasing function on $\R^+$ such
that $|D( z)| \leqslant B(y_\G(z))$. Then there exists a constant
$c_\G
>0$, depending only on $\G$ and defined in (\ref{cg}), so that
\begin{equation*}
|D(\sb z)| \leqslant B\bigl( (c_\G + 1/c_\G)( y+ 1/y) \bigr)
\end{equation*}
for any cusp $\cb$ and any $z \in \H$.
\end{lemma}

\begin{proof}
Using Lemma \ref{alem1}, which is proved in the appendix to this
paper, we have the bound
\begin{equation*}
|D(\sb z)| \leqslant B(y_\G(\sb z)) \leqslant B\bigl( (c_\G +
1/c_\G)( y+ 1/y) \bigr),
\end{equation*}
as claimed.
\end{proof}

It follows from (\ref{qmn2}), Lemma \ref{bnd71} and (\ref{yr2}) that
\begin{eqnarray*}
Q_\ca^{m,n}(\sb z,s) & \ll & \log^{m+n}\left((c_\G + 1/c_\G)( y+ 1/y) +e\right) \cdot
\left((c_\G + 1/c_\G)( y+ 1/y)\right)^\sigma    \\
& \ll & \log^{m+n}( y+ 1/y) \cdot  (y^{\sigma} + y^{-\sigma}).
\end{eqnarray*}
Substituting this bound into (\ref{qe2}) then using (\ref{fab})
and (\ref{log2}), we obtain
\begin{equation*}
E_\ca^{m,n}(\sb z,s) \ll \log^{m+n}( y+ 1/y) \cdot  (y^{\sigma} + y^{-\sigma}).
\end{equation*}
As a result, we arrive at the following conclusion.

\begin{lemma} \label{lem22} For all $z$ in $\H$, $s$ in $\C$ with  $\Re(s)>1$, and $m,n \geqslant 0$,
$$
E_\ca^{m,n}(\sb z,s) \ll y^{\sigma+\epsilon} + y^{-\sigma-\epsilon}
$$
 with an implied constant depending only on $\epsilon>0$, $\sigma$, $m$, $n$, $f$, $g$ and $\G$.
\end{lemma}

\subsection{A family of Green functions}

Set
$$
u(z,w)=\frac{|z-w|^2}{4 \,\text{Im}\,z\,\text{Im}\,w}.
$$
This function $u$ is closely related to the hyperbolic distance
between  $z$ and $w$ and satisfies the relation $u(\sigma z,\sigma
w) = u(z,w)$ for all $\sigma$ in $\SL_2(\R)$; see
\cite[(1.4)]{Iw}.
 We also need the Green function for the Laplacian:
$$
G_a(u)=\frac{1}{4\pi}\int_0^1(t(1-t))^{a-1}(t+u)^{-a}\,dt
$$
for $u,a >0$. This function is discussed in \cite[Section
1.7]{Iw}. The Green function $G_{a}(u)$ is smooth as a function of
$u$ except for a logarithmic singularity at $u=0$ which is
independent of $a$. It can be shown that, for $0<b<a$,
\begin{equation}\label{gab}
G_{ab}(u):=G_a(u)-G_b(u) \ll \frac{1}{(u+1)^b}.
\end{equation}

The next result is from \cite[Theorem 2.3]{O'S1} as well as
\cite[Section 1.9]{Iw} and is a restatement of the fact that the
resolvent of the Laplacian can be expressed as an integral
operator with kernel given by the above Green function.

\begin{theorem} \label{thmt}
If $\theta(z):{\H} \rightarrow {\C}$ is an eigenfunction of
$\Delta$ with eigenvalue $\lambda$ that satisfies $\theta(z) \ll
y^\sigma + y^{-\sigma}$ for $\sigma > 0$ then, when $a > \sigma +
1$, then
$$
\frac{-\theta(w)}{\lambda
+a(1-a)}=\int_{\H}G_a(u(w,z))\theta(z)\,d\mu(z),
$$
where  $d\mu(z)$ is the hyperbolic invariant measure $ dxdy/y^2$.
\end{theorem}

 Now we can exploit the fact that $E_{\ca}^{m,n}(z,s)$ is an eigenfunction of the
Laplacian, (\ref{del}),
along with Lemma \ref{lem22} and Theorem \ref{thmt} to write
$$
\frac{-E_{\ca}^{m,n}(z,s)}{(a(1-a)-s(1-s))}=\int_{\H}G_a(u(z,z'))E_{\ca}^{m,n}(z',s)\,d\mu
(z')
$$
for $1<$ Re$(s)<a-2$. With Lebesgue's theorems on monotone and
dominated convergence and the notation in Lemma \ref{transform},
we find
\begin{eqnarray*}
\lefteqn{\frac{-E_{\ca}^{m,n}(z,s)}{(a(1-a)-s(1-s))}
= \int_{\F}\sum_{\gamma \in \Gamma} G_a(u(z,\gamma z'))E_{\ca}^{m,n}(\g z',s)\,d\mu (z')}\\
& = & \int_{\F}\sum_{\gamma \in \Gamma} G_a(u(z,\gamma z'))
\left(\sum_{i=0}^m \sum_{j=0}^n
\binom{m}{i}\binom{n}{j}S_{m-i,n-j}(\g)E^{i,j}_{\ca}(z,s)\right)\,d\mu
(z').
\end{eqnarray*}
Therefore, for $1<$ Re$(s)< a-2$  we get
\begin{eqnarray}
\lefteqn{\frac{-E_{\ca}^{m,n}(z,s)}{(a(1-a)-s(1-s))}= \int_{\F} G_a^{0,0}(z,z')
E_{\ca}^{m,n}( z',s)\,d\mu (z')} \nonumber \hspace{35mm}\\
&+&\sum_{ 0 \leqslant i \leqslant m \atop {0 \leqslant i \leqslant
m \atop (i,j) \neq (m,n)}} \binom{m}{i}\binom{n}{j} \int_{\F}
G_a^{m-i,n-j}(z,z') E_{\ca}^{i,j}( z',s)\,d\mu (z') \label{emn21}
\end{eqnarray}
where, for $z \not\equiv z' (\text{mod}\ \Gamma)$, we have written
\begin{eqnarray*}
G_a^{i,j}(z,z') & =&  \sum_{\g} S_{i,j}(\g) G_a(u( z, \g z'))\\
&=& \sum_{\g} S_{i,j}(\g^{-1}) G_a(u(\g z, z')).
\end{eqnarray*}
The function $G_a^{i,j}(z,z')$ will inherit a logarithmic
singularity at $z \equiv z' (\text{mod}\ \Gamma)$. To eliminate
the singularity, we consider the difference
$$
G_{ab}^{i,j}(z,z'):=G_a^{i,j}(z,z')-G_b^{i,j}(z,z')= \sum_{\g} S_{i,j}(\g) G_{ab}(u( z, \g z'))
$$
for $0<b<a$.
To ease the notation we write $G_a$ and $G_{ab}$ for $G_a^{0,0}$ and $G_{ab}^{0,0}$
respectively, as in \cite[(5.1)]{Iw}. The next Proposition follows from \cite[Theorem 5.3, Lemma 5.4]{Iw}.

\begin{prop}
For all  real $a,b$ with $1<b<a$ we have that $G_{ab}(z,z') :\H
\times \H \rightarrow {\C}$ converges uniformly to a smooth
function of $z,z'$ which satisfies
$$
G_{ab}(z,z') \ll \left[ y_\G(z)  y_\G(z')\right]^{a}
$$
 where the implied constant depends only on $a,b$ and $\G$.
\end{prop}

We need a similar result for $G_{ab}^{i,j}$.

\begin{prop} \label{gij}
For all nonnegative integers $(i,j) \neq (0,0)$ and all real $a,b$
with $1<b<a$ we have that $G_{ab}^{i,j}(z,z') :\H \times \H
\rightarrow {\C}$ converges uniformly to a continuous function of
$z,z'$ which satisfies
$$
G_{ab}^{i,j}(z,z') \ll \left[ y_\F(z)  y_\F(z')\right]^{1-b+\epsilon}
$$
for any $\epsilon>0$ where the implied constant depends only on $\epsilon,a,b,i,j,f,g$ and $\G$.
\end{prop}
\begin{proof}
We begin by noting that
\begin{eqnarray*}
\frac{1}{(u( z,  z')+1)^b} &  = & \frac{ 4^b\Im (z)^b \Im( z')^b }{(|z- z'|^2 + 4  \Im(z) \Im( z'))^b} \\
& = &  \frac{ 4^b\Im (z)^b \Im( z')^b }{( \Re(z- z')^2 +   \Im(z+ z')^2 )^b}.
\end{eqnarray*}
 So we have, recalling (\ref{gab}),
\begin{eqnarray*}
G_{ab}^{i,j}( \sa z, z') & = & \sum_{ \g \in \G} S_{i,j}(\g) G_{ab}(u( \sa z, \g z'))\\
& \ll & \sum_{\g \in \G} |S_{i,j}(\g)| \frac{1}{(u( z, \sa^{-1} \g z')+1)^b}\\
& \ll & \sum_{ \g \in \G_\ca\backslash\G \atop l \in \Z}
 |S_{i,j}(\g)| \frac{ 4^b\Im (z)^b \Im(\sa^{-1}\g z')^b }{( \Re(z-\sa^{-1}\g z'+l)^2 +   \Im(z+\sa^{-1}\g z')^2 )^b}.
\end{eqnarray*}
Now since
$$
\int_{-\infty}^\infty \frac{dx}{(x^2+y^2)^b} = \sqrt{\pi} \frac{\G(b-1/2)}{\G(b) } \frac 1{y^{2b-1}}
$$
we obtain
\begin{eqnarray*}
G_{ab}^{i,j}(\sa z,z')
& \ll & \sum_{ \g \in \G_\ca\backslash\G} |S_{i,j}(\g)| \frac{\Im (z)^b \Im(\sa^{-1}\g z')^b }
{\Im(z+\sa^{-1}\g z')^{2b-1}}\\
& \ll & \Im (z)^{1-b} \sum_{ \g \in \G_\ca\backslash\G} |S_{i,j}(\g)|  \Im(\sa^{-1}\g z')^b.
\end{eqnarray*}
The proposition now follows as in the proof of Proposition \ref{prop23}.
 \end{proof}

\subsection{Applying Fredholm Theory}

We must carry out some modifications to equation (\ref{emn21}) before we can use the next result
from \cite[A.4]{Iw}.

\begin{theorem} (Fredholm) \label{fred}
Assume $\int_{\F}\,d\mu (z') =V < \infty$ and that $K(z,z')$ is
bounded and integrable on ${\F} \times {\F}$.  Then for all
$\lambda \in {\mathbb C}$ there exist $D(\lambda)$ and
$D_{\lambda}(z,z')$ entire in $\lambda$ with the following
property:  If $q(z)$ is any bounded integrable function on ${\F}$
and if $h(z)$ (defined on ${\F}$) satisfies
$$
h(z)=q(z)+\lambda\int_{\F}K(z,z')h(z')\,d\mu (z')
$$
then $h(z)$ is uniquely determined and given by the formula
$$
h(z)=q(z)+\frac{\lambda}{D(\lambda)}\int_{\F}D_{\lambda}(z,z')q(z')\,d\mu
(z')
$$
when $D(\lambda)\not= 0$.
\end{theorem}

Our situation is slightly complicated by dependence on a parameter
$s$ which will be contained in a compact set ${{\mathbf S}}\subset
{\mathbb C}$. We will set $K(z,z')=K_s(z,z')$ so $D(\lambda)$ and
$D_\lambda(z,z')$ vary with $s$. Also $q(z)=q(z,s)$ and
$\lambda=\lambda(s)$, where $q$, $\lambda$ and $K$ are analytic
functions of $s$ on ${{\mathbf S}}$. If $h(z,s)$ is analytic in
some neighborhood ${{\mathbf S}'}\subset{{\mathbf S}}$ and
satisfies
\begin{equation}\label{h31}
h(z,s)=q(z,s)+\lambda\int_{\F}K_s(z,z')h(z',s)\,d\mu (z')
\end{equation}
for all $s \in {{\mathbf S}'}$ then by Theorem \ref{fred}
$$
h(z,s)=q(z,s)+\frac{\lambda}{D(\lambda)}\int_{\F}D_{\lambda}(z,z')q(z',s)\,d\mu
(z')
$$
for all $s \in {{\mathbf S}'}$, where $D(\lambda)\not= 0$. (We
have assumed that $K_s(z,z')$ and $q(z,s)$ are uniformly bounded
in ${{\mathbf S}}$.) However, the right side of this last equation
will be meromorphic in the larger domain ${{\mathbf S}}$, which is
the method we shall use to obtain the analytic continuation of
$E^{m,n}(z,s)$.

\vskip .10in

To apply Fredholm's Theorem to the integral equation (\ref{emn21}) we need to
carry out some steps to ensure that the kernel $K$ and function $q$  are bounded.
The first is to take the difference of (\ref{emn21}) at $a$ and $b$ (for $3<b<a$) and
 eliminate the singularity at $z \equiv z'$.
Therefore, for $1<$ Re$(s)< b-2$,
\begin{equation}\label{h32}
E_{\ca}^{m,n}(z,s) \nu_{ab}(s)= q^{m,n}(z,s)+\int_{\F}
G_{ab}(z,z') E_{\ca}^{m,n}( z',s)\,d\mu (z')
\end{equation}
for
$$
q^{m,n}(z,s)= \sum_{ 0 \leqslant i \leqslant m \atop { 0 \leqslant
i \leqslant m \atop (i,j) \neq (m,n)}} \binom{m}{i}\binom{n}{j}
\int_{\F} G_{ab}^{m-i,n-j}(z,z') E_{\ca}^{i,j}( z',s)\,d\mu (z')
$$
when $\nu_{ab}(s)=(a(1-a)-s(1-s))^{-1}-(b(1-b)-s(1-s))^{-1}$.

\vskip .10in The next step is to replace the kernel $G_{ab}(z,z')$
by a kernel $G_{ab}^Y(z,z')$ that has a growing term in $z'$
removed for $z'$ in cuspidal zones (these zones depend on a large
parameter $Y$). This truncated kernel now has exponential decay as
$z'$ approaches cusps. When this new kernel is inserted into
(\ref{h32}), terms involving $\phi_{\ca \cb}^{m,n}(s)$ appear
which can be removed again by taking a certain linear combination
at $Y$, $2Y$ and $4Y$. The details of these computations are
explained in depth in \cite[Chapter 6]{Iw} and \cite[Section
3]{O'S1} so we omit the discussion. The end result of these
computations is the identity
\begin{equation}\label{h35}
h^{m,n}(z,s)=q_2^{m,n}(z,s)+\lambda(s)\int_{\F}H_s(z,z')h^{m,n}(z',s)\,d\mu
(z')
\end{equation}
when $3<\Re(s)+2 < b < a$
with
\begin{eqnarray} \nonumber
h^{m,n}(z,s) & = & (2^{s-1+a}-1)(2^{s-1+b}-1)\nu_{ab}(s)E_{\ca}^{m,n}(z,s),\\
q_2^{m,n}(z,s) & = & (2^{s-1+a}-1)(2^{s-1+b}-1) q^{m,n}(z,s),\nonumber \\
\lambda(s) & = & \frac{-1}{\nu_{ab}(s)}=\frac{(a-s)(a+s-1)(b-s)(b+s-1)}{(b-a)(a+b-1)}, \label{lam}\\
H_s(z,z') & = & \frac{ \left(G_{ab}^Y-2^{s-1}(2^a+2^b)G_{ab}^{2Y}+2^{2s-2+a+b}G_{ab}^{4Y}
\right)(z,z') }{(2^{s-1+a}-1)(2^{s-1+b}-1)}. \label{h2z}
\end{eqnarray}
The kernel $H$ still has polynomial growth at cusps in the $z$
variable. The next step is to multiply through by
$\eta(z)=e^{-\eta y_\G(z)}$ with $0<\eta<2\pi$ which controls the
growth in $z$; furthermore, the
 restriction on $\eta$ ensures that the new kernel term $\bigl(\eta(z)\eta(z')^{-1}H_s(z,z')\bigr)$
appearing below in (\ref{h36}) is bounded. Therefore, for $(m,n) \neq (0,0)$, we obtain
\begin{equation}\label{h36}
\eta(z)h^{m,n}(z,s)=\eta(z)q_2^{m,n}(z,s)+\lambda\int_{\F}
\bigl(\eta(z)\eta(z')^{-1}H_s(z,z')\bigr)
\bigl(\eta(z')h^{m,n}(z',s)\bigr)\,d\mu (z').
\end{equation}
Again, a detailed argument regarding the computations behind
(\ref{h36}) is given in  \cite{Iw} and \cite{O'S1}. Fredholm's
Theorem will be applied to (\ref{h36}) in the course of Theorem
\ref{cont2} to achieve the meromorphic continuation of $h^{m,n}$
and hence $E^{m,n}$.

\section{The continuation} \label{continuation}

\subsection{Proof of the meromorphic continuation of $E_\ca^{0,0}$}

In \cite[Sections 6.1, 6.2]{Iw} the meromorphic continuation of
$E_\ca(z,s)$ is proved.  The following proof is similar and
tailored to extend to higher orders. We also supply some details
omitted from the presentation in \cite{Iw}.

\vskip .10in
Let $\lambda$ be defined by (\ref{lam}) with fixed large constants $b<a$. Also define
$$
{\B}_{r}=\{s\in{\mathbb C}:|s| \leqslant r\}.
$$

\begin{theorem} \label{cont1}
For every ball $\B_r$ of radius $r$ and nonzero $k \in \Z$, there
exist functions $\phi_{\ca \cb}(s)$ and $\phi_{\ca \cb}(k,s)$
which are meromorphic on $\B_r$ as well as a function $A_\ca(s)$
which is analytic on $\B_r$ so that the following assertions hold
for all $s\in \B_r$:

\begin{enumerate}
\item We have the bounds
\begin{eqnarray*}
    A_\ca(s)\phi_{\ca \cb}(s) & \ll & 1,\\
    A_\ca(s)\phi_{\ca \cb}(k,s) & \ll & |k|^\sigma + |k|^{1-\sigma}.
\end{eqnarray*}

\item The Fourier expansion
\begin{equation} \label{fexp2}
E_{\ca}(\sb z,s)=\delta_{\ca \cb}
y^s+ \phi_{\ca \cb}(s)y^{1-s}+ \sum_{k\not=0}\phi_{\ca
\cb}(k,s)W_s(kz)
\end{equation}
agrees with the definition (\ref{eis}) for $\Re(s)>1$, and, for
all $z \in \H$, the series (\ref{fexp2}) converges to a
meromorphic function for $s \in\B_r$.

 \item We have the estimate
\begin{equation} \label{fbnd1}
A_\ca(s) E_{\ca}(z,s) \ll  y_\G(z)^{\sigma} + y_\G(z)^{1-\sigma}.
\end{equation}
\end{enumerate}
The implied constants in (i) and (iii) depend on $r$ and $\G$ alone.
\end{theorem}
\begin{proof}
Let
\begin{equation}\label{h00}
    h(z,s)  = \frac{(2^{s-1+a}-1)(2^{s-1+b}-1)}{2^{2s-1}-1} \nu_{ab}(s)E_{\ca}(z,s)
\end{equation}
and
\begin{equation}\label{l00}
    l(z,s) =  \frac{2^{2s-1+a-b}-1}{(2b-1)(b-s)} Y^{s-b} E_{\ca}(z,b) -
    \frac{2^{2s-1-a+b}-1}{(2a-1)(a-s)} Y^{s-a} E_{\ca}(z,a).
\end{equation}
Then it is shown in \cite[(6.8)]{Iw}, using the same method that proved (\ref{h36}), that
\begin{equation*}
\eta(z)h(z,s)=\eta(z)l(z,s)+\lambda\int_{\F}
\bigl(\eta(z)\eta(z')^{-1}H_s(z,z')\bigr)
\bigl(\eta(z')h(z',s)\bigr)\,d\mu (z')
\end{equation*}
holds for $s \in \B_{r}$ when
$$
r=b-3.
$$
(Note: (\ref{l00}) corrects a minor error in \cite[p. 84]{Iw}.) By
construction, $\eta(z)l(z,s)$ and $\eta(z)\eta(z')^{-1}H_s(z,z')$
are bounded for $(z,s) \in \H \times \B_r$ and $(z,z',s) \in \H
\times \F \times \B_r$ respectively; therefore,  it follows from
Theorem \ref{fred} that
\begin{equation}\label{hdd}
\eta(z)h(z,s)=\eta(z)l(z,s)+\frac{\lambda}{D(\lambda)}\int_{\F}
D_\lambda(z,z') \bigl(\eta(z')l(z',s)\bigr)\,d\mu (z')
\end{equation}
which is \cite[(6.10)]{Iw}. As shown in \cite[Appendix A.4]{Iw}, $D(\lambda)$
is an analytic function of $s$ and $D_\lambda(z,z')$ is analytic in $s$, piecewise
continuous in $z$, $z'$ and bounded for $(z,z',s) \in  \H \times \F  \times \B_r$.

\vskip 0.10in By combining (\ref{h00}) and (\ref{hdd}) and setting
\begin{eqnarray} \label{as}
    A_\ca(s) & = & -(2^{s-1+a}-1)(2^{s-1+b}-1)D(\lambda),\\
    A_\ca(z,s) & = & (2^{2s-1} -1) \lambda \left[D(\lambda) l(z,s) + \lambda \eta(z)^{-1}\int_{\F}
\eta(z') D_\lambda(z,z') l(z',s)\,d\mu (z') \right] \label{azs}
\end{eqnarray}
we have that $A_\ca(s)$ and $A_\ca(z,s)$ are analytic in $s$ and
$$
A_\ca(s)E_\ca(z,s) = A_\ca(z,s).
$$
It also follows from (\ref{azs}) that
\begin{equation}\label{abnd}
    A_\ca(z,s) \ll e^{\eta y_\G(z)}
\end{equation}
for $s \in \B_r$ with an implicit constant depending only on $r$ and $\G$.

\vskip 0.10in
We know the Fourier expansion (\ref{fexp2}) is valid for $\Re(s)>1$. So, for these $s$ values,
\begin{equation}\label{pke1}
    \phi_{\ca \cb}(k,s) = \frac{1}{2\sqrt{|k|y} K_{s-1/2}(2\pi |k|y)}
    \int_0^1 E_\ca(\sb(x+iy),s) e^{-2\pi ikx} \, dx.
\end{equation}
If we choose $y =1/\sqrt{|k|}$ (for later convergence results) and
replace $E_\ca(z,s)$ by $A_\ca(z,s)/A_\ca(s)$ in (\ref{pke1}) we
find
\begin{equation}\label{pke2}
    \phi_{\ca \cb}(k,s) = \frac{1}{A_\ca(s) 2|k|^{1/4}
    K_{s-1/2}(2\pi \sqrt{|k|})} \int_0^1 A_\ca(\sb(x+i/\sqrt{|k|}),s) e^{-2\pi ikx} \, dx.
\end{equation}
This yields the meromorphic continuation of the Fourier coefficients
$\phi_{\ca \cb}(k,s)$ (and similarly $\phi_{\ca \cb}(s)$) to all $s$ in $\B_r$.
But we do not yet know that (\ref{fexp2}) holds in this larger domain. To show that
it does, our task becomes proving bounds on $\phi_{\ca \cb}(k,s)$ that ensure (\ref{fexp2})
 converges to a meromorphic function on $\B_r$.

\vskip 0.10in In the arguments that follow we need additional
estimates for the $K$-Bessel functions, which we state here. In
Lemma \ref{whitaker} at the end of this section  we will prove
that
\begin{equation} \label{kb1}
2\sqrt{u} K_{s-1/2}(2\pi u) \ll e^{-2\pi u}\left( u^{r+3} + u^{-r-3} \right)
\end{equation}
for all $u>0$ and $s$ in $\B_r$ with an implied constant that
depends only on $r$. Also the asymptotic \cite[(B.36)]{Iw} shows
that there exists an absolute constant $C$ so that, for  any $s$
in $\B_r$ and $u \geqslant Cr^2$,
\begin{equation} \label{kb2}
2\sqrt{u} K_{s-1/2}(2\pi u) \geqslant  e^{-2\pi u}/2.
\end{equation}

\vskip 0.10in Continuing, we now apply Lemma \ref{bnd71} to
equation (\ref{abnd}) which produces the bound
\begin{equation}\label{abnd2}
    A_\ca(\sb(x+i/\sqrt{|k|}),s) \ll e^{\eta (c_\G + 1/c_\G)( |k|^{1/2} + |k|^{-1/2})}.
\end{equation}
 Then using (\ref{kb2}) and (\ref{abnd2}) in (\ref{pke2}) we
 arrive at the inequality
\begin{equation*}
    A_\ca(s)\phi_{\ca \cb}(k,s) \ll  e^{2 \pi \sqrt{|k|}} \cdot e^{\eta (c_\G + 1/c_\G) \sqrt{|k|}}
\end{equation*}
for $s \in \B_{r}$ and $|k| \geqslant C^2r^4$.

\vskip 0.10in Setting $y=Cr^2/|k|$ in (\ref{pke1}) instead of
$y=|k|^{-1/2}$ and using (\ref{abnd}), (\ref{kb2}) and Lemma
\ref{bnd71} shows that, for all $k$,
\begin{eqnarray} \label{pb1}
    A_\ca(s)\phi_{\ca \cb}(k,s) & \ll &  e^{2 \pi Cr^2} \int_0^1 |A_\ca(\sb(x+i Cr^2/|k|),s)| \, dx \\
    & \ll &  e^{2 \pi Cr^2}  e^{\eta y_\G(\sb(x+i Cr^2/|k|)}  \nonumber \\
    & \ll &  e^{2 \pi Cr^2} e^{\eta (c_\G+1/c_\G)\left[Cr^2/|k| + |k|/(Cr^2)\right]} \nonumber
\end{eqnarray}
for  $s \in \B_{r}$ and implied constants depending on $r$ and
$\G$ alone.  Therefore, for $|k| < C^2r^4$ we find that
$A_\ca(s)\phi_{\ca \cb}(k,s)  \ll 1$. Hence, for all $k \neq 0$,
\begin{equation}
    A_\ca(s)\phi_{\ca \cb}(k,s) \ll  e^{2 \pi D \sqrt{|k|}}  \label{abnd3}
\end{equation}
with $D=1+c_\G+ 1/c_\G$. Similar arguments show that
\begin{equation} \label{want}
    A_\ca(s)\phi_{\ca \cb}(s) \ll 1.
\end{equation}

\vskip .10in Next employ the estimate (\ref{kb1}), which is proved
in Lemma \ref{whitaker}, to see that
\begin{equation*}
    |W_s(kz)| \ll e^{-2\pi |k|y}\left( (|k|y)^{r+3} + (|k|y)^{-r-3} \right)
\end{equation*}
for $s \in \B_{r}$. As a consequence of this inequality and
(\ref{abnd3})
\begin{eqnarray}\nonumber
   \lefteqn{ A_\ca(s)\sum_{k\neq 0} \phi_{\ca \cb}(k,s)W_s(kz)  \ll
   \sum_{k=1}^\infty \left( (ky)^{r+3} + (ky)^{-r-3} \right) e^{2 \pi
   \left[D \sqrt{k}- k y \right]}}\hspace{0.8 in} \\
& \ll & e^{-2\pi y}\left[ y^{r+3} + \frac{1}{y^{r+3}}+ e^{8 \pi D^2/y}
\left( \frac{1}{y^{r-1}}+\frac{1}{y^{r+5}}\right) + e^{- \pi y}\left( y^2 +
\frac{1}{y^{r+4}}\right)\right] . \label{abnd5}
\end{eqnarray}
The last estimate (\ref{abnd5}) follows from Lemma \ref{lll} at
the end of this section. We now see that, for all $z \in \H$,
(\ref{fexp2}) converges uniformly to a meromorphic function of $s
\in \B_{r}$, proving part $(ii)$ of the theorem. Also
(\ref{abnd5}) implies that
\begin{equation}\label{abnd6}
    A_\ca(\sb z,s)= A_\ca(s)E_\ca(\sb z,s)  =  A_\ca(s)\bigl(
    \delta_{\ca \cb} y^s+ \phi_{\ca \cb}(s) y^{1-s}\bigr) +
    O(e^{-2\pi y} y^{r+3} )
    \end{equation}
as $y \to \infty$ and, consequently,
\begin{equation}
A_\ca(z,s)= A_\ca(s)E_\ca(z,s) \ll  y_\G(z)^{|\sigma-1/2| +1/2}  \label{abnd7}
\end{equation}
for all $z \in \H$. The bound (\ref{abnd7}) is a significant
improvement of (\ref{abnd}), and part $(iii)$ of the theorem
follows easily from (\ref{abnd7}). Using (\ref{abnd7}) in
(\ref{pb1}) with Lemma \ref{bnd71} produces
\begin{eqnarray*}
    A_\ca(s)\phi_{\ca \cb}(k,s) & \ll &  |k|^{|\sigma-1/2| +1/2} + |k|^{-|\sigma-1/2| -1/2} \\
    & \ll & |k|^\sigma  + |k|^{1-\sigma}.
\end{eqnarray*}
This bound, along with (\ref{want}), shows part $(i)$ of the
theorem and the proof of Theorem \ref{cont1} is complete.
\end{proof}

\begin{lemma} \label{lll} For any $D,r \geqslant 0$ and $y > 0$
$$
\sum_{k=1}^\infty k^r e^{2\pi(D \sqrt{k} - yk)} \ll  e^{-2\pi
y}\left(1 + \frac{e^{-\pi y}}{y^{r+1}} + \frac{e^{8\pi D^2/
y}}{y^{2r+2}} \right)
$$
for an implied constant depending on $r$ and $D$ alone.
\end{lemma}
\begin{proof}
We may write the left side above as
$$
e^{2\pi(D-y)}+ \sum_{k=1}^\infty (k+1)^r e^{2\pi(D\sqrt{k+1} -yk)} e^{-2\pi y}
$$
which is bounded by
\begin{equation}\label{first1}
e^{-2\pi y} \left(e^{2\pi D}+ \sum_{k=1}^\infty (2k)^r e^{2\pi(
D\sqrt{2k} -yk)} \right).
\end{equation}
The upper bound $\sqrt{2} D\sqrt{k} -yk \leqslant -yk/2$ holds if
and only if $k \geqslant 8D^2/y^2$ so
\begin{eqnarray}\nonumber
    \sum_{k=1}^\infty k^r e^{2\pi(\sqrt{2} D\sqrt{k} -y k)}
    & = & \sum_{k \leqslant 8D^2/y^2} k^r e^{2\pi(\sqrt{2}
    D\sqrt{k} -y k)} + \sum_{k > 8D^2/y^2} k^r e^{2\pi(\sqrt{2}
    D\sqrt{k} -y k)}\\ \nonumber
      & \leqslant & \sum_{k \leqslant 8D^2/y^2} (8D^2/y^2)^r e^{8\pi D^2/y} +
      \sum_{k > 8D^2/y^2} k^r e^{-\pi ky}\\
      & \leqslant & \left(\frac{8D^2}{y^2} \right)^{r+1} e^{8\pi D^2/y}  +
      \sum_{k =1}^\infty k^r e^{-\pi ky}.\label{second1}
\end{eqnarray}
Use the formula
$$
\int_0^\infty e^{-yt} t^{s-1} \, dt = \frac{\G(s)}{y^s}
$$
for $\Re(s)>0$, as in \cite[(2.4)]{CO3}, to see that
\begin{equation}\label{third1}
\sum_{k =1}^\infty k^r e^{-\pi ky} \leqslant e^{-\pi y}
\left( 1 + 2^r e^{-\pi y} +\frac{3^r \G(r+1)}{(\pi y)^{r+1}} \right)
\end{equation}
for all $r \geqslant 0$. The proof is completed by combining
inequalities (\ref{first1}), (\ref{second1}) and (\ref{third1}).
\end{proof}

\begin{lemma} \label{bessel}
For $\sigma= \Re(s) \geqslant 1$ and $y>0$,
$$
|K_{s-1/2}(y)| \leqslant \kappa(s) e^{-y}\left( y^{\sigma -1/2} + y^{-\sigma +1/2}\right)
$$
where $\kappa(s)$ is the continuous function $\kappa(s)= \pi^{1/2}
3^{\sigma-1} (1+|\G(2\sigma -1)|)/|\G(s)|$.
\end{lemma}
\begin{proof}
For $\sigma>0$ and $y>0$ we have the integral representation
$$
K_{s-1/2}(y)= \frac{\sqrt{\pi}}{\G(s)}\left( \frac y2 \right)^{s-1/2} \int_1^\infty (t^2-1)^{s-1} e^{-ty} \, dt
$$
from \cite[p. 205]{Iw}. Suppose $\sigma \geqslant 1$. It is easy to see that
$$
\left| \int_1^2 (t-1)^{s-1}(t+1)^{s-1} e^{-ty} \, dt \right| \leqslant 3^{\sigma -1} e^{-y}
$$
and
\begin{eqnarray*}
\left| \int_2^\infty (t-1)^{s-1}(t+1)^{s-1} e^{-ty} \, dt \right| & \leqslant &
e^{-y}\int_1^\infty u^{\sigma-1}(u+2)^{\sigma-1} e^{-uy} \, du \\
& \leqslant & e^{-y}3^{\sigma-1}\int_0^\infty u^{2\sigma -2} e^{-uy} \, du \\
& \leqslant & e^{-y}3^{\sigma-1} \G(2\sigma -1) y^{1-2\sigma}.
\end{eqnarray*}
Putting these estimates together finishes the proof.
\end{proof}

\begin{lemma} \label{whitaker} For all $r, y>0$ and all $s \in \B_r$
$$
W_s(y) := 2\sqrt{y} K_{s-1/2}(2\pi y) \ll e^{-2\pi y}\left( y^{r+3} + y^{-r-3} \right)
$$
with an implied constant that depends solely on $r$.
\end{lemma}
\begin{proof}
Denote by $\mathcal S_m$ the vertical strip of all $s \in \C$ with
$m \leqslant \Re(s) <m+1$.
It follows from Lemma \ref{bessel} that
\begin{eqnarray*}
|K_{s-1/2}(y)| & \leqslant & \kappa_2(s) e^{-y}\left( y^{5/2} + y^{-5/2}
\right) \ \ \text{ when }\ \ s \in \mathcal S_2,\\
|K_{s-1/2}(y)| & \leqslant & \kappa_1(s) e^{-y}\left( y^{3/2} + y^{-3/2}
\right) \ \ \text{ when }\ \ s \in \mathcal S_1.
\end{eqnarray*}
for continuous functions $\kappa_i(s)$. Then, using the identity \cite[p. 204]{Iw},
$$
K_{s-1/2}(y) = \frac{2s+1}{y} K_{s+1/2}(y) - K_{s+3/2}(y)
$$
we find
$$
|K_{s-1/2}(y)|  \leqslant  \kappa_0(s) e^{-y}\left( y^{5/2} + y^{-5/2}
\right) \ \ \text{ when }\ \ s \in \mathcal S_0
$$
and continuing we get
$$
|K_{s-1/2}(y)|  \leqslant  \kappa_{-1}(s) e^{-y}\left( y^{3/2} + y^{-7/2}
\right) \ \ \text{ when }\ \ s \in \mathcal S_{-1}
$$
and in general, for positive integers $m$,
$$
|K_{s-1/2}(y)|  \leqslant  \kappa_{-m}(s) e^{-y}\left( y^{5/2} + y^{-m-5/2}
\right) \ \ \text{ when }\ \ s \in \mathcal S_{-m}.
$$
Therefore, for all $r>0$ and $s \in \B_r$,
$$
|K_{s-1/2}(y)|  \ll   e^{-y}\left( y^{r+5/2} + y^{-r-7/2} \right)
$$
and the proof follows.
\end{proof}

\subsection{Proof of the meromorphic continuation of $E^{m,n}$}

One of the main results in this paper is the generalization of
Theorem \ref{cont1} to higher orders.

\begin{theorem} \label{cont2}
For every ball $\B_r$ of radius $r$, nonzero integer $k \in \Z$,
and all integers $m,n \geqslant 0$ there exist functions
$\phi^{m,n}_{\ca \cb}(s)$ and $\phi^{m,n}_{\ca \cb}(k,s)$ which
are meromorphic on $\B_r$, as well as a function $A_\ca(s)$ which
is analytic on $\B_r$ so that the following assertions hold for
all $s\in \B_r$:

\begin{enumerate}
\item We have the bounds
\begin{eqnarray*}
    A^{m+n+1}_\ca(s)\phi^{m,n}_{\ca \cb}(s) & \ll & 1,\\
    A^{m+n+1}_\ca(s)\phi^{m,n}_{\ca \cb}(k,s) & \ll & (\log^{m+n} |k| +1)(|k|^\sigma + |k|^{1-\sigma}).
\end{eqnarray*}

\item The Fourier expansion
\begin{equation} \label{fexpmn}
E^{m,n}_{\ca}(\sb z,s)=\delta^{m,n}_{0,0} \cdot \delta_{\ca \cb}
y^s+ \phi^{m,n}_{\ca \cb}(s)y^{1-s}+ \sum_{k\not=0}\phi^{m,n}_{\ca
\cb}(k,s)W_s(kz)
\end{equation}
 agrees with the definition (\ref{eis}) for $\Re(s)>1$ and, for all $z \in \H$, converges to a meromorphic function for  $s \in\B_r$.

 \item We have the estimate
\begin{equation} \label{fbnd}
A^{m+n+1}_\ca(s) E^{m,n}_{\ca}(z,s) \ll  y_\F(z)^{|\sigma-1/2|+1/2}.
\end{equation}
\end{enumerate}
The implied constants in $(i)$ and $(iii)$ depend on $r,m,n,f,g$ and $\G$ alone.
\end{theorem}

\begin{proof}
Our proof involves induction on $m+n$ with the first step given in
Theorem \ref{cont1}. The remainder of this proof establishes the
induction step.

\vskip 0.10in To begin, we shall now make an additional adjustment
to (\ref{h36}) in order to eliminate possible poles in $s$. Set
$q^{m,n}_A(z,s)=A_{\ca}(s)^{m+n}q^{m,n}_2(z,s)$ and
$h^{m,n}_A(z,s)=A_{\ca}(s)^{m+n}h^{m,n}(z,s)$. For $1<\Re(s)$ and
$s \in \B_r$ (with $r=b-3$ as before),
\begin{equation}\label{h37}
\eta(z)h^{m,n}_A(z,s)=\eta(z)q^{m,n}_A(z,s)+\lambda\int_{\F}
\bigl(\eta(z)\eta(z')^{-1}H_s(z,z')\bigr)
\bigl(\eta(z')h^{m,n}_A(z',s)\bigr)\,d\mu (z').
\end{equation}
We need to check that $\eta(z)q^{m,n}_A(z,s)$ is bounded before we
may apply Fredholm's Theorem (Theorem \ref{fred}) to (\ref{h37}).
Using Proposition \ref{gij} and, by induction, part $(iii)$ of the
theorem,
\begin{eqnarray*}
q^{m,n}_A(z,s) & = & A_{\ca}(s)^{m+n}(2^{s-1+a}-1)(2^{s-1+b}-1)\sum_{ 0 \leqslant i \leqslant m \atop {
0 \leqslant i \leqslant m \atop
(i,j) \neq (m,n)}}
\binom{m}{i}\binom{n}{j} \int_{\F} G_{ab}^{m-i,n-j}(z,z') E_{\ca}^{i,j}( z',s)\,d\mu (z')\\
& \ll & \sum_{ 0 \leqslant i \leqslant m \atop {
0 \leqslant i \leqslant m \atop
(i,j) \neq (m,n)}}
 \left|A_{\ca}(s)^{m+n-i-j-1}\right| \int_{\F} \left| G_{ab}^{m-i,n-j}(z,z')
 A_{\ca}(s)^{i+j+1} E_{\ca}^{i,j}( z',s) \right|\,d\mu (z')\\
& \ll & \sum_{ 0 \leqslant i \leqslant m \atop {
0 \leqslant i \leqslant m \atop
(i,j) \neq (m,n)}}
  \int_{\F} \left| y_\F(z)^{1-b+\epsilon} y_\F(z')^{1-b+\epsilon}
  \left( y_\F(z)^{\sigma} + y_\F(z)^{1-\sigma}\right) \right|\,d\mu (z')\\
& \ll & 1.
\end{eqnarray*}

In particular, the function $q^{m,n}_A(z,s)$ is  bounded for
$(z,s)\in {\F}\times{\B}_{r}$ and continuous in the $z$ variable.
We can now  apply Fredholm's theorem to (\ref{h37}) with the
result
\begin{equation}\label{h39}
h^{m,n}_A(z,s)=q^{m,n}_A(z,s)+\frac{\lambda}{D(\lambda)}\int_{\F}\eta(z)^{-1}\eta(z')D_{\lambda}
(z,z')q^{m,n}_A(z',s)\,d\mu (z')
\end{equation}
for each $s \in {\B}_{r}$ such that $D(\lambda)\not= 0$. Hence, if we set
\begin{equation}\label{h40}
A^{m,n}_\ca(z,s)=-\lambda\left(D(\lambda) q^{m,n}_A(z,s)+\lambda
\eta(z)^{-1} \int_{\F}\eta(z')D_{\lambda}
(z,z')q^{m,n}_A(z',s)\,d\mu (z')\right)
\end{equation}
we find $A^{m+n+1}_\ca(s)E^{m,n}_\ca(z,s) = A^{m,n}_\ca(z,s)$.
Also the $\eta(z)^{-1}$ term in (\ref{h40}) implies that
\begin{equation}\label{h41}
A^{m+n+1}_\ca(s)E^{m,n}_\ca(z,s) \ll e^{\eta y_\G(z)}.
\end{equation}
Now we argue as in Theorem \ref{cont1}, the only complication
being that $E^{m,n}_\ca$ is not automorphic. Since (\ref{h41}) is
also true for all smaller $m,n$ by induction, we have that
\begin{equation*}
    A^{m+n+1}_\ca(s)Q^{m,n}_\ca(z,s)  \ll  \log^{m+n}(y_\F(z)+e) \cdot e^{\eta y_\G(z)}
\end{equation*}
for $z \in \F$. Hence, for all $z \in \H$,
\begin{equation*}
    A^{m+n+1}_\ca(s)Q^{m,n}_\ca(z,s)  \ll  \log^{m+n}(y_\G(z)+e) \cdot e^{\eta y_\G(z)}
\end{equation*}
since the left side is automorphic and $\log^{m+n}(y_\F(z)+e)$ is an increasing function of $y_\F(z)$.
With Lemma \ref{bnd71},
\begin{eqnarray*}
    A^{m+n+1}_\ca(s)Q^{m,n}_\ca(\sb z,s) & \ll & \log^{m+n}((c_\G+1/c_\G)(y+1/y)+e)
    \cdot e^{\eta (c_\G+1/c_\G)(y+1/y)}\\
& \ll & \log^{m+n}(y+1/y) \cdot e^{2\pi (c_\G+1/c_\G)(y+1/y)}.
\end{eqnarray*}
Consequently, with (\ref{qe2}) and for all $z \in \H$,
\begin{eqnarray*}
    A^{m+n+1}_\ca(s)E^{m,n}_\ca(\sb z,s) & \ll & \log^{m+n}(y+1/y) \cdot e^{2\pi (c_\G+1/c_\G)(y+1/y)}\\
& \ll &  e^{2\pi D'(y+1/y)}
\end{eqnarray*}
for $D'=c_\G+1/c_\G+1$, say. Now arguing similarly to (\ref{pke1})
- (\ref{want}) we see that the Fourier coefficients of
$E^{m,n}_\ca$ may be continued and satisfy the bounds
\begin{eqnarray*}
    A^{m+n+1}_\ca(s)\phi_{\ca \cb}^{m,n}(s) & \ll & 1,\\
A^{m+n+1}_\ca(s)\phi_{\ca \cb}^{m,n}(k,s) & \ll & e^{2 \pi D'' \sqrt{|k|}}
\end{eqnarray*}
with $D''=2+c_\G+ 1/c_\G$. Thus, as in (\ref{abnd5}), we have
shown that (\ref{fexpmn}) converges to a meromorphic function for
$s \in \B_r$ and that
\begin{equation}\label{indbound}
    A^{m+n+1}_\ca(s)E^{m,n}_\ca(z,s)  \ll  y_\F(z)^{|\sigma-1/2|+1/2}.
\end{equation}
Inductively, it has been shown that the bound (\ref{indbound})
also is true for all smaller $m,n$ and therefore, again with
(\ref{qe}),
\begin{equation*}
    A^{m+n+1}_\ca(s)Q^{m,n}_\ca(z,s)  \ll  \log^{m+n}(y_\F(z)+e) \cdot y_\F(z)^{|\sigma-1/2|+1/2}.
\end{equation*}
The above result is true with $y_\F$ replaced by $y_\G$ and, employing Lemma \ref{bnd71} once more,
\begin{equation*}
    A^{m+n+1}_\ca(s)Q^{m,n}_\ca(\sb z,s)  \ll  \log^{m+n}(y+1/y) \cdot (y+1/y)^{|\sigma-1/2|+1/2}.
\end{equation*}
With (\ref{qe2}) it is now true that
\begin{equation*}
    A^{m+n+1}_\ca(s)E^{m,n}_\ca(\sb z,s)  \ll  \log^{m+n}(y+1/y) \cdot (y+1/y)^{|\sigma-1/2|+1/2}.
\end{equation*}
Use this in the analog of (\ref{pb1}) to finally see that
\begin{equation*}
    A^{m+n+1}_\ca(s)\phi_{\ca \cb}^{m,n}(k,s)  \ll  (\log^{m+n}|k| +1)  (|k|^{\sigma} + |k|^{1-\sigma}).
\end{equation*}
All parts of the theorem have been completed, establishing the induction.
\end{proof}

We have a straightforward corollary to the above results.
\begin{cor}
For every ball $\B_r$ of radius $r$ and all integers $m,n
\geqslant 0$  the following hold for all $s\in \B_r$, $z \in \H$
and $\tau \in \G$:
\begin{eqnarray}
\Delta E^{m,n}_{\ca}(z,s) & = & s(1-s)E^{m,n}_{\ca}(z,s),\label{cor1}\\
E^{m,n}_{\ca}(\tau z,s) & = & \sum_{i=0}^m\sum_{j=0}^n
\binom{m}{i}\binom{n}{j}S_{m-i,n-j}(\tau)E^{i,j}_{\ca}(z,s)\label{cor2}\\
E^{m,n}_{\ca}(\sb z,s) & = & \delta^{m,n}_{0,0} \cdot \delta_{\ca \cb}
y^s+ \phi^{m,n}_{\ca \cb}(s)y^{1-s}+ O(e^{-2\pi y}) \text{ \ \ as \ }y \to \infty.\label{cor3}
\end{eqnarray}
The implied constant in $(iii)$ will depend on $r,m,n,f,g$ and $\G$ alone.
\end{cor}
\begin{proof}
To prove (\ref{cor1}) apply the Laplacian $\Delta$ directly to the
Fourier expansion (\ref{cor3}). The relation (\ref{cor2}) clearly
follows from (\ref{emn}) by analytic continuation. To show
(\ref{cor3}) use part $(i)$ of Theorem \ref{cont2},
$$
W_s(kz) \ll e^{-2\pi ky}
$$
for $ky \geqslant r^2$ (which follows from the asymptotic \cite[(B.36)]{Iw}) and Lemma \ref{whitaker}.
\end{proof}

\section{The functional equation}

We now prove that $\mathcal E^{m,n}$ satisfies the functional
equations stated in (\ref{funeq3}) and (\ref{phieq}) which will
yield the proof of Theorem \ref{thm2}.  Recall from Theorem
\ref{cont2} that $E^{m,n}_{\ca}(z,s)$ and $\phi^{m,n}_{\ca
\cb}(s)$ are meromorphic functions of $s \in \C$.

\begin{theorem} For all $m,n \geqslant 0$ and all $s\in \C$, $z\in \H$
\begin{eqnarray}\label{funb}
 \sum_{0 \leqslant i \leqslant m \atop
0 \leqslant j \leqslant n} \binom{m}{i} \binom{n}{j}
\Phi^{i,j}(1-s) \mathcal E^{m-i,n-j}(z,s) & = & \mathcal E^{m,n}(z,1-s) \\
\sum_{0 \leqslant i \leqslant m \atop
0 \leqslant j \leqslant n} \binom{m}{i} \binom{n}{j}
\Phi^{i,j}(1-s) \Phi^{m-i,n-j}(s) & = & \begin{cases}
\id_r & \text{if \ \ $m=n=0$,} \\
0 & \text{otherwise.}
\end{cases} \label{phib}
\end{eqnarray}
\end{theorem}

\begin{proof}  First set
\begin{equation*}
\mathcal J^{m,n}(z,s):=\sum_{0 \leqslant i \leqslant m \atop
0 \leqslant j \leqslant n} \binom{m}{i} \binom{n}{j}
\Phi^{i,j}(1-s) \mathcal E^{m-i,n-j}(z,s).
\end{equation*}
When $\g \in \G$ we see by (\ref{cor2}) that
\begin{eqnarray} \nonumber
\mathcal J^{m,n}(\g z,s) & = & \sum_{0 \leqslant i \leqslant m \atop
0 \leqslant j \leqslant n} \binom{m}{i} \binom{n}{j}
\Phi^{i,j}(1-s)
\sum_{0 \leqslant c \leqslant m-i \atop
0 \leqslant d \leqslant n-j} \binom{m-i}{c} \binom{n-j}{d} S_{m-i-c,n-j-d}(\g)
\mathcal E^{c,d}(z,s)\\ \nonumber
 & = & \sum_{0 \leqslant c \leqslant m \atop
0 \leqslant d \leqslant n} \binom{m}{c} \binom{n}{d} S_{m-c,n-d}(\g) \sum_{0 \leqslant i \leqslant c \atop
0 \leqslant j \leqslant d} \binom{c}{i} \binom{d}{j}
\Phi^{i,j}(1-s) \mathcal E^{c-i,d-j}(z,s)\\
& = & \sum_{0 \leqslant c \leqslant m \atop
0 \leqslant d \leqslant n} \binom{m}{c} \binom{n}{d} S_{m-c,n-d}(\g)  \mathcal J^{c,d}(z,s). \label{fund}
\end{eqnarray}
To prove the theorem we use induction on $m+n$ with the base case,
$(m,n)=(0,0)$, being the well-known functional equation for the
classical non-holomorphic Eisenstein series. 
For the induction step we now assume that
\begin{equation}\label{func}
\mathcal E^{c,d}(z,1-s) = \mathcal J^{c,d}(z,s)
\end{equation}
for all $c,d$ with $c+d<N$ and suppose that $m+n=N$. 
Combine (\ref{fund}) and (\ref{func}) to see that
\begin{eqnarray*}
\mathcal J^{m,n}(\g z,s) & = & \mathcal J^{m,n}(z,s) + \sum_{0 \leqslant c \leqslant m \atop
{0 \leqslant d \leqslant n \atop (c,d) \neq (m,n) }} \binom{m}{c} \binom{n}{d} S_{m-c,n-d}(\g)  \mathcal J^{c,d}(z,s)\\
& = & \mathcal J^{m,n}(z,s) + \mathcal E^{m,n}(\g z,1-s) - \mathcal E^{m,n}(z,1-s).
\end{eqnarray*}
Therefore the function $\mathcal J^{m,n}(z,s)  - \mathcal
E^{m,n}(z,1-s)$ is $\G$-invariant in the $z$ variable. We also see
from (\ref{cor1}) that $\mathcal J^{m,n}(z,s)$ is an eigenfunction
on the Laplacian with eigenvalue $s(1-s)$ and, with (\ref{fbnd}),
must be bounded by $y_\G(z)^{|\sigma-1/2|+1/2}$. Therefore, by a
result of Selberg (see \cite[Lemma 6.4]{Iw}), $\mathcal
J^{m,n}(z,s)- \mathcal
E^{m,n}(z,1-s)$ is equal to a linear combination of the classical
($m=n=0$) non-holomorphic Eisenstein series:
\begin{equation} \label{jep}
 \mathcal J^{m,n}(z,s)  - \mathcal E^{m,n}(z,1-s) = \Psi^{m,n}(s)\mathcal E^{0,0}(z,s)
\end{equation}
for some $r \times r$ matrix $\Psi^{m,n}(s)$. Considering the Fourier expansions of both sides of (\ref{jep}) and  comparing the
coefficients of $y^s$  shows that
$$
\Psi^{m,n}(s)=0
$$
and proves (\ref{funb}). Comparing the coefficients of $y^{1-s}$
then yields (\ref{phib}).  We have now completed the induction
step and the proof is complete.
\end{proof}

\section{Further examples}

Consider a group $\Gamma$ with $r$ cusps.  Taking $(m,n) = (2,1)$,
the action by $\Gamma$ as described in Theorem \ref{thm1} becomes
the identity
$$
\left( \begin{matrix} \mathcal E^{2,1} \\ \mathcal E^{2,0} \\ \mathcal E^{1,1} \\ \mathcal E^{1,0} \\
\mathcal E^{0,1} \\ \mathcal E^{0,0} \end{matrix} \right)(\g z,s) = \overline{\chi}(\g)\left( \begin{matrix}
1 & -\overline{\s{\g}{g}} & -2 \s{\g}{f} & 2 \s{\g}{f} \overline{\s{\g}{g}} & \s{\g}{f}^2 &  -\s{\g}{f}^2
\overline{\s{\g}{g}}\\
 & 1 & & -2\s{\g}{f}  & & \s{\g}{f}^2 \\
& & 1 & -\overline{\s{\g}{g}} & -\s{\g}{f}  & \s{\g}{f} \overline{\s{\g}{g}} \\
& & & 1 & & -\s{\g}{f} \\
& & & & 1 & -\overline{\s{\g}{g}}  \\
& & & & & 1
 \end{matrix} \right)
\left( \begin{matrix} \mathcal E^{2,1} \\ \mathcal E^{2,0} \\ \mathcal E^{1,1} \\ \mathcal E^{1,0} \\ \mathcal
E^{0,1} \\ \mathcal E^{0,0} \end{matrix} \right)(z,s)
$$
where, as in (\ref{eg1}), each matrix entry above is multiplied by the $r \times r$ identity matrix $\id_r$.
The functional equation in Theorem \ref{thm2} is the relation
$$
\left( \begin{matrix} \mathcal E^{2,1} \\ \mathcal E^{2,0} \\ \mathcal E^{1,1} \\ \mathcal E^{1,0} \\
\mathcal E^{0,1} \\ \mathcal E^{0,0} \end{matrix} \right)(z,1-s) = \left( \begin{matrix} \Phi^{0,0}
& \Phi^{0,1}& 2\Phi^{1,0}& 2\Phi^{1,1}& \Phi^{2,0}& \Phi^{2,1} \\
 & \Phi^{0,0} & & 2\Phi^{1,0} & & \Phi^{2,0} \\
& & \Phi^{0,0} &\Phi^{0,1} &\Phi^{1,0} &\Phi^{1,1} \\
& & & \Phi^{0,0} & &\Phi^{1,0}  \\
& & & & \Phi^{0,0} & \Phi^{0,1}  \\
& & & & & \Phi^{0,0}
 \end{matrix} \right)(1-s)
\left( \begin{matrix} \mathcal E^{2,1} \\ \mathcal E^{2,0} \\ \mathcal E^{1,1} \\ \mathcal E^{1,0} \\ \mathcal
E^{0,1} \\ \mathcal E^{0,0} \end{matrix} \right)(z,s).
$$
As throughout the paper, each entry $\mathcal E$ is a column
vector of size $r \times 1$, and each entry $\Phi$ is an $r \times
r$ matrix.

\vskip .10in Consider the specific example when $\G=\G_0(p)$ is
the Hecke congruence group of prime level $p$.  The group
$\G_{0}(p)$ has two inequivalent cusps which my be uniformized to
be at $\infty$ and $0$.  In this case, we have that for each
$(m,n)$, the vector $\mathcal E^{m,n}$ is simply the order pair
${}^t(E_\ci ^{m,n}, E_0 ^{m,n})$.  Taking $(m,n)=(2,0)$, and
omitting $z$ from the notation, Theorem \ref{thm2} yields the
identity
\begin{equation}\label{Hecke}
\left( \begin{matrix}   E_\ci^{2,0} \\  E_0^{2,0} \\ E_\ci^{1,0} \\
E_0^{1,0} \\ E_\ci^{0,0} \\  E_0^{0,0} \\   \end{matrix}
\right)(1-s) = \left( \begin{matrix}   \Phi_{\ci \ci}^{0,0} & \Phi_{\ci
0}^{0,0}& 2\Phi_{\ci \ci}^{1,0} & 2\Phi_{\ci 0}^{1,0}&
\Phi_{\ci \ci}^{2,0} & \Phi_{\ci 0}^{2,0} \\
\Phi_{0 \ci}^{0,0} & \Phi_{0 0}^{0,0}& 2\Phi_{0 \ci}^{1,0} & 2\Phi_{0 0}^{1,0}& \Phi_{0 \ci}^{2,0}
& \Phi_{0 0}^{2,0} \\
 & &\Phi_{\ci \ci}^{0,0} & \Phi_{\ci 0}^{0,0}& \Phi_{\ci \ci}^{1,0} & \Phi_{\ci 0}^{1,0}
 \\
 & &\Phi_{0 \ci}^{0,0} & \Phi_{0 0}^{0,0}& \Phi_{0 \ci}^{1,0} & \Phi_{0 0}^{1,0}\\
  & & & & \Phi_{\ci \ci}^{0,0} & \Phi_{\ci 0}^{0,0} \\
 & & & & \Phi_{0 \ci}^{0,0} & \Phi_{0 0}^{0,0}
\end{matrix} \right)(1-s)
\left( \begin{matrix}   E_\ci^{2,0} \\  E_0^{2,0} \\ E_\ci^{1,0} \\
E_0^{1,0} \\ E_\ci^{0,0} \\  E_0^{0,0} \\  \end{matrix} \right)
(s).
\end{equation}
As  in \cite{O'S1}, each entry in the scattering matrix in
(\ref{Hecke}) can be expressed as a type of Kloosterman sum
twisted by modular symbols. It certainly would be interesting to
find explicit formulas for the entries of (\ref{Hecke}).  The
matrix identity (\ref{funeq2}) implies a number of relations for
these scattering matrix entries.

\begin{appendix}

\section{Invariant height} \label{app}

For any Fuchsian group $\G$ we define a constant $c_\G>0$ as follows. Let
$$
\mathcal C_{\ca \cb} = \left\{ |c| : (\smallmatrix * & * \\ c & *
\endsmallmatrix) \in \sa^{-1} \G \sb, \ c \neq 0\right\}.
$$
As described in \cite[Sections 2.5, 2.6]{Iw} and \cite[Lemma
1.25]{Sh}, there is a smallest element in the set $\mathcal C_{\ca
\cb}$ which we shall denote by $c(\ca,\cb)$. With this, we let
\begin{equation} \label{cg}
c_\G:= \max_{\ca, \cb} \{ c(\ca,\cb)^{-2}\}.
\end{equation}
An important point, shown in \cite[Section 2.4]{Iw}, is that for any $\g \in \G$
\begin{equation} \label{c0}
\sa^{-1} \g \sb=(\smallmatrix * & * \\ 0 & * \endsmallmatrix) \implies
\ca=\cb \ \ \text{ and } \ \ \sa^{-1} \g \sb=(\smallmatrix 1 & l \\ 0 & l \endsmallmatrix)
\end{equation}
for some $l \in \R$. The next result is similar to the second part of \cite[Lemma 1.25]{Sh}.

\begin{lemma} \label{alem1} For all $z\in \H$ and every cusp $\cb$
\begin{equation*}
y_\G(\sb z) \leqslant (c_\G+ 1/c_\G)(y+1/y).
\end{equation*}
\end{lemma}

\begin{proof}
We first prove the following fact about the invariant height
function $y_\G(w)$: If $y_\G(w)>Y$ then, for any cusp $\cb$,
either $\Im(\sb^{-1} w)>Y$ or $\Im(\sb^{-1} w)< c_\G/Y$. To see
why this is true note that $y_\G(w)>Y$ implies that there exists
$\ca$ and $\g \in \G$ so that $\Im(\sa^{-1}\g w)>Y$. Hence
\begin{equation} \label{yY}
|j(\sa^{-1} \g \sb, z)|^2 <\frac yY
\end{equation}
for $z=\sb^{-1}w$. Let $\g'=\sa^{-1} \g \sb=(\smallmatrix * & * \\
c & d \endsmallmatrix)$. Now if $c=0$ then (\ref{c0}) implies that
$\g'=(\smallmatrix 1 & * \\ 0 & 1 \endsmallmatrix)$ and therefore
$\Im(\sb^{-1} w)>Y$. Otherwise, if $c \neq 0$, we have from
(\ref{yY}) that
\begin{equation*}
|cy|^2 <\frac yY
\end{equation*}
and $y<1/(|c|^2Y)$. That implies $y<c_\G /Y$ as required.

\vskip .10in The proof of the lemma is now straightforward. The
contrapositive of the first result is that if there exists a cusp
$\cb$ so that
\begin{equation*}
c_\G/Y \leqslant \Im(\sb^{-1} w) \leqslant Y
\end{equation*}
then $y_\G(w) \leqslant Y$. In other words, for any cusp $\cb$,
\begin{equation} \label{cyY}
c_\G/Y \leqslant y \leqslant Y \implies y_\G(\sb z) \leqslant Y.
\end{equation}
Choosing $Y=(c_\G+ 1/c_\G)(y+1/y)$ in (\ref{cyY}) finishes the proof.
\end{proof}

For completeness we also prove a lower bound for the invariant height of any point $z$ in $\H$.

\begin{lemma} \label{alem2} For every $z\in \H$
\begin{equation*}
y_\G(z) \geqslant  \begin{cases}
c_\G & \text{if \ \ $c_\G \leqslant 1/2$,} \\
\sqrt{c_\G -1/4} & \text{if \ \ $c_\G > 1/2$.}
\end{cases}
\end{equation*}
\end{lemma}
\begin{proof}
Suppose  $y_\G(w)<Y$ then, for all cusps $\ca$, $\cb$ and all $\g
\in \G$ we have  $\Im(\sa^{-1}\g \sb z )<Y$ for $z=\sb^{-1}w$.
Hence
\begin{equation} \label{yY2}
|j(\sa^{-1} \g \sb, z)|^2 >\frac yY.
\end{equation}
We need to make a careful choice of $\ca$, $\cb$ and $\g$ to get
the desired bounds. First choose $\ca$, $\cb$ and $\g_0$ so that
$\sa^{-1} \g_0 \sb=(\smallmatrix * & * \\ c & d \endsmallmatrix)$
with $c_\G = 1/|c|^2$. Recall that $\G_\cb$, the subgroup of $\G$
consisting of elements that fix $\cb$, is generated by a primitive
element $\g_\cb$ that satisfies $\sb^{-1} \g_\cb^l \sb =
(\smallmatrix 1 & l \\ 0 & 1 \endsmallmatrix)$. Hence
$$
\sa^{-1} \g_0 \g_\cb^l \sb=(\smallmatrix * & * \\ c & d +cl \endsmallmatrix).
$$
Choose $l$ so that $|x+d/c +l| \leqslant 1/2$ and put $\g = \g_0 \g_\cb^l$ in (\ref{yY2}) to get
\begin{equation*}
|j(\sa^{-1} \g \sb, z)|^2 = |c|^2\left(|x+d/c+l|^2 +y^2\right)>\frac yY
\end{equation*}
and therefore
\begin{equation} \label{yY4}
1/4 > \frac {c_\G y}Y -y^2.
\end{equation}
Complete the square in (\ref{yY4}) and replace $Y$ by $c_\G/ \sqrt{u^2+1}$ for a new parameter $u$ to get
\begin{equation*}
y=\Im(\sb^{-1} w)  > (u+\sqrt{u^2+1})/2.
\end{equation*}
Since $y_\G(w) \geqslant \Im(\sb^{-1} w)$, our work to this point has shown that
\begin{equation} \label{yY6}
y_\G(w)<\frac{c_\G}{\sqrt{u^2+1}} \implies y_\G(w)> (u+\sqrt{u^2+1})/2
\end{equation}
for any $u \geqslant 0$.

\vskip 0.10in If $c_\G \leqslant 1/2$ then taking $u=0$ in
(\ref{yY6}) shows that $y_\G(w)<c_\G$ is impossible. If $c_\G >
1/2$ then
\begin{equation*}
\frac{c_\G}{\sqrt{u^2+1}} = (u+\sqrt{u^2+1})/2
\end{equation*}
has a positive solution, namely $u=(2c_\G -1)/\sqrt{4c_\G -1}$.
Substituting this value into (\ref{yY6}), we arrive at the bound
\begin{equation*}
y_\G(w)<\sqrt{c_\G -1/4} \implies y_\G(w)> \sqrt{c_\G -1/4}
\end{equation*}
and deduce that $y_\G(w)<\sqrt{c_\G -1/4}$ is impossible. This completes the proof.
\end{proof}

Let $y_\G$ be the minimum value of $y_\G(z)$. For example, when
$\G = \SL_2(\Z)$ we have $c_\G=1$ and therefore by Lemma
\ref{alem2} we know $y_\G \geqslant \sqrt{3}/2$. In fact, by
examining the Ford fundamental domain for $\SL_2(\Z)$ we see that
$y_\G = \sqrt{3}/2$.

\vskip 0.10in
We conclude with the result quoted in the proof of Proposition \ref{prop23}.

\begin{lemma} \label{eminus} For $\Re(s)=\sigma >1$ and an implied constant depending on $\sigma$, $\G$ alone,
\begin{equation*}
E_\ca(z,s) - \Im (\sa^{-1} z)^s  \ll y_\F(z)^{1-\sigma}.
\end{equation*}
\end{lemma}
\begin{proof}
We know from \cite[Corollary 3.5]{Iw} that
\begin{equation*}
E_\ca(\sb z,s) = \delta_{\ca \cb} y^s + \phi_{\ca \cb}(s) y^{1-s} + O(e^{-2 \pi y})
\end{equation*}
as $y \to \infty$. Hence
\begin{equation*}
E_\ca(\sa z,s) - \Im (\sa^{-1} \sa z)^s  \ll y^{1-\sigma}
\end{equation*}
as $y \to \infty$ and, for $\ca \neq \cb$
\begin{eqnarray*}
E_\ca( \sb z,s) - \Im (\sa^{-1} \sb z)^s & = & - \Im (\sa^{-1} \sb z)^s + \phi_{\ca \cb}(s) y^{1-s} + O(e^{-2 \pi y}) \\
& \ll & \Im (\sa^{-1} \sb z)^\sigma + y^{1-\sigma}.
\end{eqnarray*}
By (\ref{c0}) we have $\sa^{-1} \sb=(\smallmatrix * & * \\ c & d \endsmallmatrix)$ with $c \neq 0$ and
$$
\Im (\sa^{-1} \sb z)^\sigma = \frac{y^\sigma}{|cz+d|^{2\sigma}} = \frac{y^\sigma}{(|cx+d|^2 + |cy|^2)^{\sigma}}
\leqslant \frac{y^\sigma}{|cy|^{2\sigma}} \leqslant (c_\G)^\sigma y^{-\sigma}.
$$
The result follows.
\end{proof}

\end{appendix}

\bibliography{refdata}

\vskip .25in \noindent Jay Jorgenson \newline Department of
Mathematics \newline City College of New York
\newline Convent Avenue at 138th Street \newline New York, NY
10031 \newline e-mail: jjorgenson@mindspring.com

\vskip .15in \noindent Cormac O'Sullivan \newline Department of
Mathematics \newline Bronx Community College \newline University
Avenue and West 181st Street \newline Bronx, NY 10453 \newline
e-mail: cormac12@juno.com

\end{document}